\documentclass[a4paper, 11pt]{amsart}
\usepackage[utf8]{inputenc}

\usepackage[english]{babel}

\usepackage{fancyhdr}
\usepackage{amsmath,amssymb,amsthm}
\usepackage{mathrsfs}
\usepackage{tikz}

\usepackage{a4wide}
\usepackage[text={17.5cm,24cm}]{geometry}
\usepackage{graphicx}
\usepackage[arrow, matrix, curve]{xy}
\usepackage{MnSymbol}

\numberwithin{equation}{section}

\theoremstyle{definition}
\newtheorem{df}{Definition}[section]
\newtheorem{rem}[df]{Remark}
\newtheorem{ex}[df]{Example}

\theoremstyle{plain}
\newtheorem{lemma}[df]{Lemma}
\newtheorem{prop}[df]{Proposition}
\newtheorem{thm}[df]{Theorem}
\newtheorem{cor}[df]{Corollary}

\newenvironment{bew}{\begin{proof}[Proof]}{\end{proof}}

\renewcommand{\a}{\mathbb A}
\renewcommand{\c}{\mathbb C}
\newcommand{\f}{\mathbb F}
\newcommand{\ka}{\mathbb K}

\newcommand{\p}{\mathbb P}

\renewcommand{\t}{\mathbb T}
\newcommand{\z}{\mathbb Z}

\newcommand{\F}{\mathcal F}

\newcommand{\M}{\mathfrak M}

\newcommand{\Jon}{\text{\textbf{Jon}}}
\newcommand{\Aff}{\text{\textbf{Aff}}}
\newcommand{\chr}{\text{char}}
\newcommand{\reg}{\text{reg}}

\DeclareMathOperator{\Aut}{Aut}

\DeclareMathOperator{\supp}{\text{supp}}

\DeclareMathOperator{\id}{\text{id}}
\DeclareMathOperator{\ML}{\text{ML}}

\DeclareRobustCommand{\ie}{i.\,e.~}

\newcommand{\nlin}{\unitlength1mm\begin{picture}(0,9.25)
                       \put(0,0.75){\line(0,1){8.5}}
                      \end{picture}}

\newcommand{\vlin}[1]{\hspace{0.75mm}\unitlength1mm\begin{picture}(#1,0)
                       \put(0,0){\line(1,0){#1}}
                      \end{picture}\hspace{0.75mm}\rule[-3mm]{0mm}{4mm}}

\newcommand{\lin}{\vlin{8.5}}

\newcommand{\lllin}{\vlin{15}}

\newcommand{\co}[1]{\unitlength1mm\begin{picture}(0,8)
    \put(0,0){\circle{1.5}}
    \put(0,3){\makebox(0,5)[b]{$#1$}}
                      \end{picture}}

\newcommand{\mybox}{\unitlength1mm\begin{picture}(0,1.5)
    \put(-0.75,-0.75){\line(0,1){1.5}}
    \put(-0.75,-0.75){\line(1,0){1.5}}
    \put(0.75,0.75){\line(0,-1){1.5}}
    \put(0.75,0.75){\line(-1,0){1.5}}
    \end{picture}}

\newcommand{\xbox}{\unitlength1mm\begin{picture}(0,1.5)
    \put(0,0){$\mybox$}
    \put(-0.75,0){\line(1,0){1.5}}
    \put(0,-0.75){\line(0,1){1.5}}
    \end{picture}}
    
\newcommand{\xboxo}[1]{\unitlength1mm\begin{picture}(0,8)
    \put(0,0){\xbox}
    \put(0,3){\makebox(0,5)[b]{$#1$}}
                      \end{picture}}

\newcommand{\cu}[1]{\unitlength1mm\begin{picture}(0,8)
    \put(0,0){\circle{1.5}}
    \put(0,-7){\makebox(0,4)[b]{$#1$}}
    \end{picture}
      \rule[-7mm]{0mm}{7mm}}

\newcommand{\cou}[2]{\unitlength1mm\begin{picture}(0,8)
    \put(0,0){\circle{1.5}}
    \put(0,3){\makebox(0,5)[b]{$#1$}}
    \put(0,-8.5){\makebox(0,4)[t]{$#2$}}
      \end{picture}
      \rule[-7mm]{0mm}{7mm}}

\newcommand{\cshiftup}[2]{\unitlength1mm\begin{picture}(0,9.25)
                       \put(0,10){\cou{#1}{#2}}
                      \end{picture}}

\newcommand{\xbou}[2]{\unitlength1mm\begin{picture}(0,8)
    \put(0,0){\xbox}
    \put(0,2){\makebox(0,5)[b]{$#1$}}
   \put(0,-7){\makebox(0,5)[b]{$#2$}}
      \end{picture}
      \rule[-7mm]{0mm}{7mm}}

\newcommand{\xbshiftup}[2]{\unitlength1mm\begin{picture}(0,9.25)
                       \put(0,10){\xbou{#1}{#2}}
                      \end{picture}}

\begin{document}

\begin{center}
\LARGE{Smooth Non-Homogeneous Gizatullin Surfaces}

\parskip 20pt

\normalsize Sergei Kovalenko
\end{center}

\vspace{10pt}

\begin{sloppypar}
\noindent\small ABSTRACT. Quasi-homogeneous surfaces, or Gizatullin surfaces, are normal affine surfaces such that there exists an open orbit of the automorphism group with a finite complement. If the action of the automorphism group is transitive, the surface is called homogeneous. Examples of non-homogeneous Gizatullin surfaces were constructed in \cite{Ko}, but on more restricted conditions. We show that a similar result holds under less constrained assumptions. Moreover, we exhibit examples of smooth affine surfaces with a non-transitive action of the automorphism group whereas the automorphism group is huge. This means that it is not generated by a countable set of algebraic subgroups and that its quotient by the (normal) subgroup generated by all algebraic subgroups contains a free group over an uncountable set of generators.
\end{sloppypar}

\tableofcontents

\begin{sloppypar} 

\section{Introduction}

Let $\ka$ be an algebraically closed ground field $\ka$. \emph{Quasi-homogeneous surfaces} or \emph{Gizatullin surfaces} were studied by Danilov and Gizatullin (\cite{DG1}, \cite{DG2} and \cite{DG3}). These are normal affine surfaces over $\ka$ which, except for $\ka^* \times \ka^*$, satisfy one of the equivalent conditions in the following theorem:\\

\begin{thm}\label{TheoremGizatullin} (see \cite{Gi} and \cite{Du} for the normal case) For a normal affine surface that is non-isomorphic to $\ka^* \times \ka^*$, the following conditions are equivalent:
\begin{itemize}
\item[(1)] The automorphism group $\Aut(V)$ acts on $V$ with an open orbit $O$, such that the complement $V \backslash O$ is finite ($O$ is called the \textit{big orbit} of $\Aut(V)$).
\item[(2)] $V$ admits a smooth compactification by a smooth zigzag $D$. In other words, $V = X \backslash D$, where $X$ is a complete surface smooth along $D$ and $D$ is a linear chain of smooth rational curves with simple normal crossings.
\end{itemize} 
\end{thm}

Recall that the Makar-Limanov invariant $\ML(V)$ of an affine surface $V$ is defined to be the intersection of all kernels of locally nilpotent derivations of the coordinate ring $\ka[V]$. Assuming now that $\chr(\ka) = 0$, these two conditions are, except for $V = \ka^* \times \a^1$, equivalent to $\ML(V)$ being trivial, that is, $\ML(V) = \ka$. Normal affine surfaces $V$ satisfying one of the equivalent conditions of Theorem \ref{TheoremGizatullin} are called \emph{quasi-homogeneous surfaces} or \emph{Gizatullin surfaces}. Moreover, $V$ is called \emph{homogeneous} if $O$ coincides with $V_\reg$, that is, if $\Aut(V)$ acts transitively on $V_\reg$. In particular, the automorphism group of a Gizatullin surface $V$ is quite large compared to surfaces in general.

In positive characteristic, examples of quasi-homogeneous surfaces which are not homogeneous were early known, see \cite{DG1}. Gizatullin formulated in \cite{Gi} his conjecture that every smooth quasi-homogeneous surface is already homogeneous if $\ka$ has characteristic $0$. In \cite{Ko} the author constructed counterexamples to this conjecture.
 
The aim of this article is to strengthen the main result of \cite{Ko} and to construct more general families of quasi-homogeneous surfaces, which are not homogeneous. This provides a criterion for a quasi-homogeneous surface to be non-homogeneous under less constrained assumptions. In particular, we determine finite subsets that are invariant under the action of $\Aut(V)$.

Let $(X, D)$ be an SNC-completion of a Gizatullin surface $V$ so that $V = X \backslash D$ and $D$ is a simple normal crossing divisor. By applying suitable birational transformations we can transform $D$ into \emph{standard form} (\cite{DG2}). The latter means that $D = C_0 \cup \cdots \cup C_n$ is a chain of smooth rational curves with either $C_0^2 = C_1^2 = 0$ and $C_i^2 \leq -2$ for $i \geq 2$ if $n \geq 4$ or with $C_i^2 = 0$ for all $i$ if $n \leq 3$. The sequence $[[C^2_0, C^2_2, \dots, C^2_n]]$ is called the \emph{type} of $D$. Up to reversion, the standard form of the boundary divisor $D$ is an invariant of the abstract isomorphism type of $V$ (\cite{FKZ1C} Cor. 3.33'). However, since this invariant provides little information about the surface in general, it is more convenient to consider a stronger invariant, the so called \emph{extended divisor} $D_{\text{ext}}$, which is defined as follows. Since $C^2_0 = C^2_1 = 0$, we obtain two $\p^1$-fibrations $\Phi_0 := \Phi_{|C_0|}: \tilde{X} \to \p^1$ and $\Phi_1 := \Phi_{|C_1|}: \tilde{X} \to \p^1$, where $\tilde{X}$ is the minimal resolution of singularities of $X$. By \cite{FKZ2}, Lemma 2.19, $\Phi_0$ has at most one degenerate fiber, without lost of generality the fiber over $0$, and the \emph{extended divisor} of $(X, D)$ is

$$D_{\text{ext}} := C_0 \cup C_1 \cup \Phi_0^{-1}(0).$$

By construction, the extended divisor $D_{\text{ext}}$ always contains the boundary divisor $D$, and is a tree (\cite{FKZ3}, Prop. 1.11). The connected components of $D_{\text{ext}} - D$ are called \emph{feathers}. We denote them by $F_{i, j}$, $2 \leq i \leq n$, $j \in \{1, \dots, r_i \}$ and assume that $F_{i, j}$ is attached to the curve $C_i$ at the point $P_{i, j}$. Moreover, if $X$ is smooth (hence $\tilde{X} = X$), the feathers are irreducible. Furthermore, the \emph{Matching Principle} (cf. \cite{FKZ4}) provides a natural bijection between feathers $F_{i, j}$ of $(X, D)$ and feathers $F^\vee_{i^\vee, j}$ of the completion $(X^\vee, D^\vee)$, which is obtained by reversing the boundary zigzag.

We concentrate on smooth Gizatullin surfaces which admit a \emph{$(-1)$-completion} (see Def. \ref{DfMinusOneCompletion}). These are, by definition, Gizatullin surfaces which admit a standard completion $(X, D)$ such that every feather of $D_{\text{ext}}$ has self-intersection number $-1$. The following theorem (see Theorem \ref{MainTheorem}), which is the main result of this article, provides a wide class of smooth non-homogeneous Gizatullin surfaces:

\begin{thm}\label{MainTheoremIntroduction} Let $V$ be a smooth Gizatullin surface which admits a $(-1)$-completion $(X, D)$. Let $A_i = \{ P_{i, 1}, \dots, P_{i, r_i} \} \subseteq C_i \backslash (C_{i - 1} \cup C_{i + 1}) \cong \c^*$ be the base point set of the feathers $F_{i, j}$. For a finite subset $A \subseteq \c^*$, we denote by $G(A)$ the group $\{ \alpha \in \c^* \mid \alpha \cdot A = A \}$, and for an inner boundary component $C_i$ of $D$, such that $i \not\in \mathfrak{E}_D \cup \mathfrak{E}^\vee_{D^\vee}$ (for the definition of the subsets $\mathfrak{E}_D, \mathfrak{E}^\vee_{D^\vee} \subseteq \{ 2, \dots, n \}$ see Def. \ref{DfExceptionalComponent}), we let $B_{i, 1}, \dots, B_{i, m_i}$ be the orbits of the $G(A_i)$-action on $A_i$. If the configuration invariant $Q(X, D) = (Q_2, \dots, Q_n)$ of $V$ (see \ref{SectionMatchingPrinciple}) is not symmetric, we let 

\begin{equation}
O_{i, j} := \bigcup_{1 \leq l \leq r_i; P_{i, l} \in B_{i, j}} (F_{i, l} \cap F^\vee_{i, l}) \subseteq V, \quad j = 1, \dots, m_i.
\end{equation}

\noindent Otherwise, we let for $i \leq \lfloor \frac{n}{2} \rfloor + 1$

\begin{equation}
O_{i, j} := \left( \bigcup_{1 \leq l \leq r_i; P_{i, l} \in B_{i, j}} (F_{i, l} \cap F^\vee_{i, l}) \right) \cup \left( \bigcup_{1 \leq l \leq r_i; P_{i^\vee, l} \in B_{i^\vee, j}} (F_{i^\vee, l} \cap F^\vee_{i^\vee, l}) \right) \subseteq V, \quad j = 1, \dots, m_i,
\end{equation}

\noindent where we identify $B_{i, j}$ with $B_{i^\vee, j}$ under a suitable isomorphism $C_i \backslash (C_{i - 1} \cup C_{i + 1}) \stackrel{\sim}{\to} C_{n + 2 - i} \backslash (C_{n + 3 - i} \cup C_{n + 1 - i})$. Moreover, in both cases, we let  

$$O_0 := V \backslash \left( \bigcup_{i, j} O_{i, j} \right).$$

\noindent Then the following hold:
\begin{itemize}
\item[(1)] The subsets $O_0$ and $O_{i, j}$ are invariant under the action of $\Aut(V)$. Moreover, $O_0$ contains the big orbit $O$.
\item[(2)] Let $F$ be a feather of $D_{\text{ext}}$ which is attached to $C_i$, such that either $C_i$ is an outer component or $i \in \mathfrak{E}_D \cup \mathfrak{E}^\vee_{D^\vee}$ holds. Then $F \backslash D$ is contained in $O$. 
\item[(3)] Assume, that $r_i > 0$ holds for a unique $i \in \{ 2, \dots, n \}$. Then the subsets $O_0$ and $O_{i, j}$ form the orbit decomposition of the natural action of $\Aut(V)$ on $V$. In particular, $O = O_0$ holds.
\end{itemize}
\end{thm}

We will also exhibit among these non-homogeneous Gizatullin surfaces examples of surfaces which admit a huge group of automorphisms. Here we say that an automorphism group $\Aut(V)$ is \emph{huge}, if: \\
\noindent (1) The (normal) subgroup $\Aut(V)_{\text{alg}}$ generated by all algebraic subgroups is not generated by a countable set of algebraic subgroups, and\\
\noindent (2) The quotient $\Aut(V)/\Aut(V)_{\text{alg}}$ contains a free group over an uncountable set of generators.\\
Other examples of Gizatullin surfaces with huge automorphism group were constructed in \cite{BD2}. More precisely, it is shown that if $V$ is a smooth Gizatullin surface of type $[[0, 0, -a, -b]]$ with $a, b \geq 3$, then $\Aut(V)$ is huge, if the feathers of $D_{\text{ext}}$ are attached to general points of the corresponding components. However, according to Theorem \ref{MainTheorem}, such surfaces are already homogeneous.\\

This article is structured as follows. In section 2 we recall the basic notions and tools concerning $\a^1$-fibrations and $\a^1$-fibered surfaces. We recall the notion of the extended divisor and the Matching Principle.

In section 3 we apply these tools to give a description of the correspondence fibration for a smooth Gizatullin surface in the general case. This description allows us to deduce Theorem \ref{MainTheorem}.

Finally, section 4 deals with the structure of the automorphism group of some special smooth Gizatullin surfaces. In particular, we exhibit the promised examples of non-homogeneous Gizatullin surfaces with a huge automorphism group.

\bigskip
\noindent \textbf{Acknowledgements} The author would like to thank Jérémy Blanc for helpful discussions on the structure of the automorphism groups. Furthermore, the author thanks Mathias Leuenberger for inspiring discussions on the problem of holomorphic automorphisms of Gizatullin surfaces.

\section{Preliminaries}

\bigskip
\subsection{$\a^1$-fibered surfaces and Gizatullin surfaces}

In this section we recall some basic facts about $\a^1$-fibered surfaces and, in particular, about Gizatullin surfaces. We work over the field $\ka = \c$ of complex numbers, but all results stated in this section are also valid for an arbitrary algebraically closed field of characteristic $0$. Let us recall the notion of an oriented zig (see \cite{DG2}).

\begin{df} A \emph{zigzag} $D$ on a normal projective surface $X$ is an SNC-divisor supported in the smooth locus $X_{\text{reg}}$ of $X$, with irreducible components isomorphic to $\p^1$ and whose dual graph is a chain. If $\supp(D) = \bigcup_{i = 0}^n C_i$ is the decomposition into irreducible components, one can order the $C_i$ such that

$$C_i.C_j = \begin{cases} 1, & |i - j| = 1 \\ 0, & |i - j| > 1 .\end{cases}$$

A zigzag with such an ordering is called \emph{oriented} and the sequence $[[(C_0)^2, \dots, (C_n)^2]]$ is called the \emph{type} of $D$. The same zigzag with the reverse ordering is denoted by ${}^tD$, \ie ${}^tD$ is of type $[[(C_n)^2, \dots, (C_0)^2]]$.

An \emph{oriented sub-zigzag} of an oriented zigzag is an SNC-divisor $D'$ with $\supp(D') \subseteq \supp(D)$ which is a zigzag for the induced ordering.

We say that an oriented zigzag $D$ is composed of sub-zigzags $Z_1, \dots, Z_s$, and following \cite{BD1} we denote $D = Z_1 \triangleright \cdots \triangleright Z_s$, if the $Z_i$, $1\leq i\leq s$, are oriented sub-zigzags of $D$ whose union is $D$ and the components of $Z_i$ precede those of $Z_j$ for $i < j$.
\end{df}

Surfaces completable by a zigzag were first studied by Danilov and Gizatullin (\cite{Gi}, \cite{DG2} and \cite{DG3}).

\begin{df} Normal affine surfaces $V$ satisfying one of the conditions in Theorem \ref{TheoremGizatullin} are called \emph{quasi-homogeneous surfaces} or \emph{Gizatullin surfaces}. If the big orbit $O$ of the natural action of $\Aut(V)$ on $V$ coincides with $V_\reg$, then $V$ is called \emph{homgeneous}.
\end{df}

\noindent For the rest of this article we fix the following notation:\\

\noindent \textbf{Notation:} If $V$ is a Gizatullin surface and $(X, D)$ is a completion of $V$ by a zigzag $D$, then $D = C_0 + \cdots + C_n$ and $C_i$ and $C_j$ have a non-empty intersection only for $|i - j| = 1$. In particular, the natural number $n$ \emph{always} denotes the length of the boundary zigzag $D$.

\bigskip
Let $V$ be a Gizatullin surface and $(X, D)$ be a completion of $V$ by a zigzag. We can associate a linear weighted graph $\Gamma_D$ to $(X, D)$ as follows. The vertices $v_i$, $0 \leq i \leq n$, are the boundary components $C_i$ and the weights are the corresponding self-intersection numbers $w_i := C_i^2$. In other words, $\Gamma_D$ has the form

$$\Gamma_D: \quad \cou{C_0}{w_0} \lin \cou{C_1}{w_1} \lin \ \cdots \ \lin \cou{C_n}{w_n} \quad .$$

%For a better systematic understanding of Gizatullin surfaces we introduce elementary transformations of weighted graphs.

%\begin{df} Given an at most linear vertex $v$ of a weighted graph $\Gamma$ with weight $0$ one can perform the following transformations. If $v$ is linear with neighbors $v_1, v_2$ then we blow up the edge connecting $v$ and $v_1$ in $\Gamma$ and blow down the proper transform of $v$:

%\begin{equation}\label{ElementaryTransformation1}
%\dots \ \cou{v_1}{w_1 - 1} \lin \cou{v'}{0} \lin \cou{v_2}{w_2 + 1} \ \dots \ \dasharrow \ \dots \ \cou{v_1}{w_1 - 1} \lin \cou{v'}{-1} \lin \cou{v}{-1} \lin \cou{v_2}{w_2} \ \dots \ \to \ \dots \ \cou{v_1}{w_1} \lin \cou{v}{0} \lin \cou{v_2}{w_2}
%\end{equation}

%\noindent Similarly, if $v$ is an end vertex of $\Gamma$ connected to the vertex $v_1$ then one proceeds as follows:

%\begin{equation}\label{ElementaryTransformation2}
%\dots \ \cou{v_1}{w_1 - 1} \lin \cou{v'}{0} \ \dasharrow \ \dots \ \cou{v_1}{w_1 - 1} \lin \cou{v'}{-1} \lin \cou{v}{-1} \ \to \ \dots \ \cou{v_1}{w_1} \lin \cou{v}{0}
%\end{equation}

%\noindent These operations (\ref{ElementaryTransformation1}) and (\ref{ElementaryTransformation2}) and their inverses are called \emph{elementary transformations} of $\Gamma$. If such an elementary transformation involves only an inner blow-up then we call it \emph{inner}. Thus (\ref{ElementaryTransformation1}) and (\ref{ElementaryTransformation2}) are inner whereas the inverse of (\ref{ElementaryTransformation2}) is not as it involves an outer blow-up.
%\end{df}

Applying a suitable sequence of blow-ups and blow-downs we can transform the dual graph $\Gamma_D$ of $D$ into \textit{standard form}, \ie we can achieve that $C_0^2 = C_1^2 = 0$ and $C_i^2 \leq -2$ for all $i \geq 2$ if $n \geq 4$ or $C_i^2= 0$ for all $i$ if $n \leq 3$ (see \cite{DG2}, \cite{Da}, \cite{FKZ1}). Moreover, this representation is unique up to reversion meaning that for two standard forms $[[0, 0, w_2, \dots, w_n]]$ and $[[0, 0, w'_2, \dots, w'_n]]$ either $w_i = w'_i$ or $w_i = w'_{n + 2 - i}$ holds (\cite{FKZ1C}, Cor. 3.33').

The reversion process can be described as follows. We start with a boundary divisor of type $[[0, 0, w_2, \dots, w_n]]$. Blowing up $X$ in $C_0 \cap C_1$ and contracting the proper transform of $C_1$ yields a boundary divisor of the type $[[-1, 0, w_2 + 1, \dots, w_n]]$. Repeating this procedure $|w_2|$ times we arrive at a zigzag of type $[[w_2, 0, 0, w_3, \dots, w_n]]$. In this way we can move the pair of zeros to the right and obtain finally a zigzag of type $[[w_2, \dots, w_n, 0, 0]]$. Note, that all birational transformations are centred in the boundary, \ie these transformations yield isomorphisms on the affine parts.

\begin{df}\label{DefinitionZigzag} A zigzag $D$ on a normal projective surface $X$ is called \emph{$m$-standard} (or in \emph{$m$-standard form}), if it is of type $[[0, -m, w_2, \dots, w_n]]$ with $n \geq 1$ and $w_i \leq -2$ (in the case of $n = 1$ there are no weights $w_i$).\\
\noindent An \emph{$m$-standard pair} is a pair $(X, D)$ consisting of a normal projective surface $X$ and an $m$-standard zigzag $D$ on $X$. If $m = 0$, then $(X, D)$ is called a \emph{standard pair}. A \emph{birational map} $\varphi: (X, D) \dasharrow (X', D')$ between $m$-standard pairs is a birational map $\varphi: X \dasharrow X'$ which restricts to an isomorphism $\varphi\vert_{X \backslash D}: X \backslash D \stackrel{\sim}{\to} X' \backslash D'$.
\end{df}

Let $(X, D)$ be an $m$-standard pair and let $(\tilde{X}, D) \to (X, D)$ be a minimal resolution of singularities. Since $\tilde{X}$ is rational and $C_0$ is a $0$-curve, the linear system $|C_0|$ defines a $\p^1$-fibration $\bar{\pi} = \Phi_{|C_0|}: \tilde{X} \to \p^1$. In particular, if $m = 0$, there are even two $\p^1$-fibrations $\Phi_0 := \Phi_{|C_0|}, \Phi_1 := \Phi_{|C_1|} : \tilde{X} \to \p^1$, and thus a morphism 

$$\Phi := \Phi_0 \times \Phi_1: \tilde{X} \to Q = \p^1 \times \p^1,$$

\noindent which is birational (\cite{FKZ2}, Lemma 2.19). Choosing suitable coordinates on the quadric $Q$ we can assume that $C_0 = \Phi_1^{-1}(\infty)$, $\Phi(C_1) = \{ \infty \} \times \p^1$ and $C_2 \cup \cdots \cup C_n \subseteq \Phi_1^{-1}(0)$. The divisor $D_{\text{ext}} := C_0 \cup C_1 \cup \Phi_1^{-1}(0)$ is called the \textit{extended divisor}. We also denote the full fiber $\Phi_0^{-1}(0)$ by $D_{(e)}$. For determining the structure of the extended divisor, we recall the notion of a \emph{feather}:

\begin{df} (\cite{FKZ2}, Def. 5.5)
\begin{itemize}
\item[(1)] A \emph{feather} is a linear chain 

$$F: \quad \co{B}{} \lin \co{F_1} \lin \dots \lin \co{F_s}$$

\noindent of smooth rational curves such that $B^2 \leq -1$ and $F_i^2 \leq -2$ for all $i \geq 1$. The curve $B$ is called the \emph{bridge curve}.
\item[(2)] A collection of feathers $\{ F_\rho \}$ consists of feathers $F_\rho$, $1 \leq \rho \leq r$, which are pairwise disjoint. Such a collection will be denoted by a plus box 

$$\xboxo{ \{ F_\rho \} } \quad .$$

\item[(3)] Let $D = C_0 + \cdots + C_n$ be a zigzag. A collection $\{ F_\rho \}$ is \emph{attached to a curve} $C_i$ if the bridge curves $B_\rho$ meet $C_i$ in pairwise distinct points and all the feathers $F_\rho$ are disjoint with the curves $C_j$ for $j \neq i$.
\end{itemize} 
\end{df}

\begin{lemma}\label{StructureExtendedDivisor} (\cite{FKZ3}, Prop. 1.11) Let $(\tilde{X}, D)$ be a minimal SNC completion of the minimal resolution of singularities of a Gizatullin surface $V$. Furthermore, let $D = C_0 + \cdots + C_n$ be the boundary divisor in standard form. Then the extended divisor $D_{\text{ext}}$ has the dual graph

\bigskip
$$D_{\text{ext}}: \quad \cou{0}{C_0} \lin \cou{0}{C_1} \lin \cu{C_2} \nlin \xbshiftup{ \{ F_{2, j} \} }{} \lin \dots \lin \cu{C_i} \nlin \xbshiftup{ \{ F_{i, j} \} }{} \lin \dots \lin \cu{C_n} \nlin \xbshiftup{ \{ F_{n, j} \} }{} \quad ,$$

\noindent where $\{ F_{i, j} \}$, $j \in \{1, \dots, r_i \}$, are feathers attached to the curve $C_i$. Moreover, $\tilde{X}$ is obtained from $\p^1 \times \p^1$ by a sequence of blow-ups with centers in the images of the components $C_i$, $i \geq 2$.
\end{lemma}

\noindent \textbf{Notation:} For the rest of this article, the numbers $r_2, \dots, r_n$ always denote the number of feathers attached to $C_2, \dots, C_n$.

\begin{rem} We consider the feathers $F_{i, j} := B_{i, j} + F_{i, j, 1} + \cdots + F_{i, j, k_{i, j}}$ mentioned in Lemma \ref{StructureExtendedDivisor}. The collection of linear chains $R_{i, j} := F_{i, j, 1} + \cdots + F_{i, j, k_{i, j}}$ corresponds to the minimal resolution of singularities of $V$. Thus, if $(X, D)$ is a standard completion of $V$ and $(\tilde{X}, D)$ is the minimal resolution of singularities of $(X, D)$, the chain $R_{i, j}$ contracts via $\mu: (\tilde{X}, D) \to (X, D)$ to a singular point of $V$, which is a cyclic quotient singularity. In partcular, $V$ has at most cyclic quotient singularities (see \cite{Mi}, \S 3, Lemma 1.4.4 (1) and \cite{FKZ3}, Remark 1.12).

It follows that $V$ is smooth if and only if every $R_{i, j}$ is empty, \ie if every feather $F_{i, j}$ is irreducible and reduces to a single bridge curve $B_{i, j}$ (\cite{FKZ3}, 1.8, 1.9 and Remark 1.12).
\end{rem}

In the following we abbreviate the subdivisor $\sum_{k \geq i} C_k + \sum_{j_k; k \geq i} F_{k, j_k}$ by $D^{\geq i}_{\text{ext}}$ and the subdivisor $\sum_{k > i} C_k + \sum_{j_k; k \geq i} F_{k, j_k} = D^{\geq i}_{\text{ext}} \ominus C_i$ by $D^{> i}_{\text{ext}}$.

\bigskip
A similar statement as in Lemma \ref{StructureExtendedDivisor} holds for minimal resolutions of singularities of $1$-standard completions of Gizatullin surfaces. They arise as blow-ups of the Hirzebruch surface $\f_1$.

\begin{lemma}\label{FactorizationHirzebruch}(\cite{BD1}, Lemma 1.0.7) Let $(X, D)$ be a $1$-standard pair and let $\mu: \tilde{X} \to X$ be the minimal resolution of singularities of $X$. Then there exists a birational morphism $\eta: \tilde{X} \to \f_1$, unique up to an automorphism of $\f_1$, that restricts to an isomorphism outside the degenerate fibers of $\bar{\pi} \circ \mu$, and satisfies the commutative diagram

$$\begin{xy}
  \xymatrix{
   & \tilde{X} \ar[ld]_\mu \ar[rd]^\eta \ar[dd]^{\mu \circ \bar{\pi}} & \\
  X \ar[rd]_{\bar{\pi}} &  & \f_1 \ar[ld]^\rho \\
   & \p^1 & .
  }
\end{xy}$$

\noindent Moreover, if $(X', D')$ is another $1$-standard pair with associated morphism $\eta': \tilde{X}' \to \f_1$, then $(X, D)$ and $(X', D')$ are isomorphic if and only if there exists an automorphism of $\f_1$ isomorphically mapping $\eta(\mu^{-1}_*(C_0))$ onto $\eta'(\mu'^{-1}_*(C'_0))$ and isomorphically sending the base-points of $\eta^{-1}$ (including infinitely near ones) onto those of $\eta'^{-1}$.
\end{lemma}

\bigskip
The study of automorphisms of quasi-homogeneous surfaces leads in a natural way to the study of birational maps betweens their completions. Indeed, every automorphism of a quasi-homogeneous surface $V$ can be extended to a birational automorphism of a standard completion $(X, D)$ of $V$. However, birational automorphisms can be controlled much better, since they admit decompositions into "elementary" birational maps (see Prop. \ref{FactorizationOfBirationalMaps}).

In the following we give a short description of birational maps between standard pairs as well as between $1$-standard pairs. It follows from \cite{BD1}, Lemma 2.1.1 that every birational map $\varphi: (X, D) \dasharrow (X', D')$ between $1$-standard pairs, which is not an isomorphism, has a unique base point $p \in C_0$. This base point is called the \textit{center} of $\varphi$. In general, this yields qualitatively different maps depending on whether $p \in C_0 \cap C_1$ or $p \in C_0 \backslash C_1$. We recall these two types of birational maps in the following definition:

\begin{df} Let $\varphi: (X, D) \dasharrow (X', D')$ be a birational map between $1$-standard pairs and let $D = C_0 \triangleright \cdots \triangleright C_n$ and $D' = C'_0 \triangleright \cdots \triangleright C'_n$ be the oriented boundary divisors.
\item[(1)] (Fibered modification) $\varphi$ is called a \emph{fibered map} if it restricts to an isomorphism of $\a^1$-fibered quasi-projective surfaces

$$\begin{xy}
\xymatrix{
V = X \backslash D \ar[r]^\sim_\varphi \ar[d]_{\bar{\pi}\mid_{V}} & V' = X' \backslash D' \ar[d]^{\bar{\pi}'\mid_{V'}} \\
\a^1 \ar[r]^\sim & \a^1.
}
\end{xy}$$

\noindent $\varphi$ is called \emph{fibered modification} if it is not an isomorphism.
\item[(2)] (Reversion) $\varphi$ is called \emph{reversion} if it admits a resolution of the form

$$\begin{xy}
\xymatrix{
 & (Z, \tilde{D} = C_n \triangleright \cdots \triangleright C_1 \triangleright H \triangleright C'_1 \triangleright \cdots \triangleright C'_{n'}) \ar[ld]_\sigma \ar[rd]^{\sigma'} & \\
(X, {}^tD) \ar@{-->}[rr]^\varphi &  & (X', D'),
}
\end{xy}$$

\noindent where $H$ is a zigzag with boundaries $C_0$ (left) and $C'_0$ (right) and where $\sigma: Z \to X$ and $\sigma': Z \to X'$ are smooth contractions of the sub-zigzags $H \triangleright C'_1 \triangleright \cdots \triangleright C'_{n'}$ and $C_n \triangleright \cdots \triangleright C_1 \triangleright H$ of $\tilde{D}$ onto $C_0$ and $C'_0$ respectively.
\end{df}

\begin{rem} In a similar way we define fibered modifications for $m$-standard pairs: a birational map $\varphi: (X, D) \dasharrow (X', D')$ between $m$-standard pairs is called a \emph{fibered modification} if it restricts to an isomorphism of $\a^1$-fibered quasi-projective surfaces

$$\begin{xy}
\xymatrix{
V = X \backslash D \ar[r]^\sim_\varphi \ar[d]_{\bar{\pi}\mid_{V}} & V' = X' \backslash D' \ar[d]^{\bar{\pi}'\mid_{V'}} \\
\a^1 \ar[r]^\sim & \a^1.
}
\end{xy}$$

\noindent and is not an isomorphism.
\end{rem}

By \cite{BD1}, Lemma 2.4.1, every fibered modification $(X, D) \dasharrow (X', D')$ between $1$-standard pairs is centred in $p = C_0 \cap C_1$ and every reversion $(X, D) \dasharrow (X', D')$ between $1$-standard pairs is centred in $p \in C_0 \backslash C_1$. Furthermore, the center $p$ gives the full control over the reversion:

\begin{prop}\label{UniquenessReversions} (Uniqueness of reversions, see \cite{BD1}, Prop. 2.3.7) For every $1$-standard pair $(X, D)$ and every point $p \in C_0 \backslash C_1$ there exist a $1$-standard pair $(X', D')$ and a reversion $\varphi: (X, D) \dasharrow (X', D')$, unique up to an isomorphism at the target, having $p$ as a unique proper base point. Moreover, if $\Gamma_D = [[0, -1, w_2, \dots, w_n]]$, then $\Gamma_{D'} = [[0, -1, w_n, \dots, w_2]]$.
\end{prop}

Every birational map between standard completions decomposes into "elementary maps", fibered modifications and reversions, as the following proposition states:

\begin{prop}\label{FactorizationOfBirationalMapsStandard} (\cite{Ko}, Cor. 3.3) Let $\varphi: (X, D) \dasharrow (X', D')$ be a birational map between standard pairs. Then there exists a decomposition

$$\varphi = \varphi_m \circ \cdots \circ \varphi_1: (X, D) = (X_0, D_0) \stackrel{\varphi_1}{\to} (X_1, D_1) \stackrel{\varphi_2}{\to} \cdots \stackrel{\varphi_m}{\to} (X_m, D_m) = (X', D')$$

\noindent such that each $\varphi_i$ is either a reversion or a fibered modification.
\end{prop}

\noindent Prop. \ref{FactorizationOfBirationalMapsStandard} holds as well for $1$-standard pairs. We will need the following statement in Section \ref{SectionHugeAut}, in particular, its uniqueness part:

\begin{prop}\label{FactorizationOfBirationalMaps}(\cite{BD1}, Theorem 3.0.2) Let $\varphi: (X, D) \dasharrow (X', D')$ be a birational map between $1$-standard pairs restricting to an isomorphism $X \backslash D \stackrel{\sim}{\to} X' \backslash D'$. If $\varphi$ is not an isomorphism, then it can be decomposed into a finite sequence

$$\varphi = \varphi_n \circ \cdots \circ \varphi_1: (X, D) = (X_0, D_0) \stackrel{\varphi_1}{\to} (X_1, D_1) \stackrel{\varphi_2}{\to} \cdots \stackrel{\varphi_n}{\to} (X_n, D_n) = (X', D')$$

\noindent of fibered modifications and reversions between $1$-standard pairs $(X_i, D_i)$. Moreover, such a factorization of minimal length is unique, meaning, if 

$$\varphi = \varphi'_n \circ \cdots \circ \varphi'_1: (X, D) = (X'_0, D'_0) \stackrel{\varphi'_1}{\to} (X'_1, D'_1) \stackrel{\varphi'_2}{\to} \cdots \stackrel{\varphi'_n}{\to} (X'_n, D'_n) = (X', D')$$

\noindent is another factorization of minimal length, then there exist isomorphisms of $1$-standard pairs $\alpha_i: (X_i, D_i) \to (X'_i, D'_i)$, such that $\alpha_i \circ \varphi_i = \varphi'_i \circ \alpha_{i - 1}$ for $i = 2, \dots, n$.
\end{prop}

\begin{rem} Prop. 2.3.3 in \cite{FKZ4} asserts a similar factorization of birational maps between semi-standard completions of a Gizatullin surface. Any semi-standard completion of a Gizatullin surface can be obtained from another one by a finite sequence of so called \emph{generalized reversions}, that is, reversions $(X, D) \dasharrow (X', D')$ such that $C_1^2 = 0$ and ${C'_1}^2 = 0$ do not necessarily hold. However, to pass from generalized reversions to reversions between standard pairs, one needs fibered modifications. 
\end{rem}

\bigskip
\subsection{The Matching Principle}\label{SectionMatchingPrinciple}

In the following we give a short overview over the Matching Principle for Gizatullin surfaces (see \cite{FKZ4}, section 3). Let $(X, D)$ be a standard completion of a smooth Gizatullin surface $V$, $(X, D) \dasharrow (X^\vee, D^\vee)$ be the reversion and let $D = C_0 \cup \cdots \cup C_n$ and $D^\vee = C^\vee_0 \cup \cdots \cup C^\vee_n$. Furthermore, we let $\Gamma_D = [[0, 0, w_2, \dots, w_n]]$ and we denote the corresponding extended divisors by $D_{\text{ext}}$ and $D^\vee_{\text{ext}}$, respectively. Performing inner elementary transformations on the boundary divisor $D$ we can move the pair of zeros to the right by several places. Let us abbreviate for a given integer $2 \leq t \leq n$

$$t^\vee := n + 2 - t.$$

So for every $t$ with $2 \leq t \leq n + 1$ we obtain a new completion $(W, E)$ of $V$ with boundary divisor $[[w_2, \dots, w_{t - 1}, 0, 0, w_t, \dots, w_n]]$, \ie

$$E = C^\vee_n \cup \cdots \cup C^\vee_{t^\vee} \cup C_{t - 1} \cup C_t \cup \cdots \cup C_n,$$

\noindent if we identify $C_i \subseteq X$ and $C^\vee_j \subseteq X^\vee$ with their proper transforms in $W$. In particular, we can write $E$ as $E = D^{\geq t - 1} \cup D^{\vee \geq t^\vee}$ with new weights $C^2_{t - 1} = C^{\vee 2}_{t^\vee} = 0$. Moreover, we have natural isomorphisms

\begin{eqnarray*} 
W \backslash D^{\vee \geq t^\vee} &=& W \backslash(C^\vee_n \cup \cdots \cup C^\vee_{t^\vee}) \cong X \backslash (C_0 \cup \cdots \cup C_{t - 2}),\\
W \backslash D^{\geq t - 1} &=& W \backslash(C_{t - 1} \cup \cdots \cup C_n) \cong X^\vee \backslash (C^\vee_0 \cup \cdots \cup C^\vee_{t^\vee - 1}).\\
\end{eqnarray*}

\begin{df}(\cite{FKZ4}, Def. 3.3.3) The map

$$\psi := \Phi_{|C_{t - 1}|}: W \to \p^1$$

\noindent is called the \emph{correspondence fibration} for the pair $(C_t, C^\vee_{t^\vee})$.
\end{df}

The Matching Principle provides a natural correspondence between feathers of $D_{\text{ext}}$ and those of $D^\vee_{\text{ext}}$. Recall that for a given feather $F$ of $D_{\text{ext}}$, a boundary component $C_\mu$ is called \emph{mother component of} $F$ if the feather $F$ is created by a blowup on $C_\mu$ during the blowup process $X \to Q$ (see \cite{FKZ3}, 2.3).

\begin{prop}\label{MatchingFeathers} (cf. \cite{FKZ4}, Lemma 3.3.4, Cor. 3.3.5 and Lemma 3.3.6) Let $F$ be a feather of $D_{\text{ext}}$ attached to the component $C_i$. Then there exists a unique feather $F^\vee$ of $D^\vee_{\text{ext}}$ which intersects $F$ in $V$ and which is attached to a component $C^\vee_j$, such that $i + j \geq n + 2$. Moreover, $F$ and $F^\vee$ intersect transversally and in a single point. If $C_\tau$ is the mother component of $F$, then $C^\vee_{\tau^\vee}$ is the mother component of $F^\vee$. 
\end{prop}

\begin{df} Feathers $F$ and $F^\vee$, which satisfy the conditions of Proposition \ref{MatchingFeathers}, are called \emph{matching feathers}.
\end{df}

\begin{rem} Note, that the condition $i + j \geq n + 2$ is essential. Indeed, every feather $F_{t - 1, \rho}$ is a section of $\psi$ and therefore it meets every fiber of $\psi$. Since it cannot intersect $D^{\geq t}_{\text{ext}}$, it meets every feather $G_{t^\vee, \sigma}$ of $D^\vee_{\text{ext}}$ with $(G_{t^\vee, \sigma})^2 = -1$ on $V$. Moreover, if the condition $i + j \geq n + 2$ does not hold, the feathers $F$ and $F^\vee$ may intersect in more than one point. For example, consider a Gizatullin surface $V$ with extended divisor

\bigskip
$$D_{\text{ext}}: \quad \cu{0} \lin \cu{0} \lin \cu{-2} \nlin \cshiftup{-1}{} \lin \cu{-3} \nlin \cshiftup{-1}{} \lin \dots \lin \cu{-3} \nlin \cshiftup{-1}{} \lin \cu{-2} \nlin \cshiftup{-1}{} \quad ,$$

\bigskip
\noindent and denote the feather attached to $C_i$ by $F_i$. Then, using the algorithm in section \ref{PresentationGizatullinSurface} below (or \cite{FKZ4}, 5.1.1), it is easy to see that $F_2 \cap F^\vee_i$ consists of $i - 1$ points, if the $F_i$ are attached to general points of $C_i \backslash (C_{i - 1} \cup C_{i + 1})$.
\end{rem}

\bigskip
\noindent \textbf{Configuration spaces and the configuration invariant.} We consider a smooth Gizatullin surface $V$ with a standard completion $(X, D)$. The sequence of weights $[[w_2, \dots, w_n]]$ (up to reversion) of the boundary divisor $D$ is a discrete invariant of the abstract isomorphism type of $V$ (\cite{FKZ1C}, Cor. 3.33'). However, two Gizatullin surfaces may be non-isomorphic, even if the dual graphs their extended divisors coincide. The reason is, that the configurations of the base points $p_{i, j} = F_{i, j} \cap C_i$  of the feathers may differ. A partial solution to this problem is a stronger continuous invariant of $V$, the \textit{configuration invariant}, which takes such configurations into account. In the following we recall the notion of the configuration invariant (see \cite{FKZ4}, Section 3).

For a natural number $s \geq 1$ we denote the configuration space of all $s$-points subsets $\{ \lambda_1, \dots, \lambda_s \} \subseteq \a^1$ by $\mathcal{M}^+_s$. We can identify $\mathcal{M}^+_s$ in a natural way with the Zariski open subset of $\a^s$:

$$\mathcal{M}^+_s \cong \a^s \backslash \{ \text{discr}(P) = 0 \}, \quad \text{where} \quad P = \prod_{j = 1}^s (X - \lambda_j),$$

\noindent see \cite{FKZ4}, 3.1.1. The group $\Aut(\a^1)$ acts on $\mathcal{M}^+_s$ in a natural way. We let

$$\M^+_s := \mathcal{M}^+_s/\Aut(\a^1).$$

\noindent Thus, $\M^+_s$ is an $(s - 2)$-dimensional affine variety.

Now, let $\mathcal{M}^*_s$ be the configuration space of all $s$-points subsets $\{ \lambda_1, \dots, \lambda_s \} \subseteq \c^* = \a^1 \backslash \{ 0 \}$. Similarly, the group $\c^*$ acts on $\mathcal{M}^*_s$ and we let 

$$\M^*_s := \mathcal{M}^*_s/\c^*.$$

Before introducing the configuration invariant we have to distinguish two types of boundary components.

\begin{df}\label{StarComponent}
\begin{itemize}
\item[(1)] For a natural number $i \in \{ 2, \dots, n \}$ $s_i$ shall denote the number of feathers of $D_{\text{ext}}$ whose mother component is $C_i$.
\item[(2)] The component $C_i$ is called a \emph{$*$-component} or \emph{inner component} if 
\begin{itemize}
\item[(i)] $D^{\geq i + 1}_{\text{ext}}$ is not contractible and
\item[(ii)] $D^{\geq i + 1}_{\text{ext}} - F_{j, k}$ is not contractible for every feather $F_{j, k}$ of $D^{\geq i + 1}_{\text{ext}}$ with mother component $C_\tau$, where $\tau < i$.
\end{itemize}
\noindent Otherwise $C_i$ is called a \emph{$+$-component} or \emph{outer component}.
\end{itemize}
\end{df}

For example, $C_2$ and $C_n$ are always $+$-components. In the following we let $\tau_i = *$ in the first case and $\tau_i = +$ in the second one.

It is not hard to see that in the blow-up process $\tilde{X} \to \p^1 \times \p^1$ ($\tilde{X}$ is a standard completion of the minimal resolution of singularities $V'$ of $V$) every $*$-component $C_i$, $3 \leq i \leq n - 1$, appears as a result of an inner blow-up of the previous zigzag, while an outer blow-up of a zigzag creates a $+$-component.

\begin{lemma}\label{StarComponentAfterReversion}(\cite{FKZ4}, Lemma 3.3.10) $C_t$ is a $*$-component if and only if $C^\vee_{t^\vee}$ is a $*$-component.
\end{lemma}

Now we are able to construct the so-called \textit{configuration invariant of} $V$. First, let $C_i$ be a $+$-component. For every feather $F_{i, j}$ with self-intersection $-1$ we let $p_{i, j}$ be its intersection point with $C_i$. Moreover, if there exists a feather $F_{k, j}$ with mother component $C_i$ and $k > i$, then we also add the intersection point $c_{i + 1} := C_i \cap C_{i + 1}$ to our collection. Note, that such a feather is unique, if it exists. Thus, the collection of points

$$p_{i, j} \in C_i, \quad 1 \leq j \leq s_i$$

\noindent is just the collection of locations on $C_i$ in which the feathers with mother component $C_i$ are born by a blow-up. These points are called \textit{base points} of the associated feathers. The collection $(p_{ij})_{1 \leq j \leq s_i}$ defines a point $Q_i$ in $\M^+_{s_i}$.

Let now $C_i$ be a $*$-component. In the same way as above we consider $Q_i$ as a collection of points on $C_i \backslash (C_{i - 1} \cup C_{i + 1})$. Note that the intersection point $c_{i + 1}$ of $C_i$ and $C_{i + 1}$ cannot belong to this collection due to Definition \ref{StarComponent} (2) (ii). Identifying $C_i \backslash (C_{i - 1} \cup C_{i + 1})$ with $\c^*$ in a way that $c_{i + 1}$ corresponds to $0$ and $c_i$ to $\infty$ we obtain a point $Q_i$ in the configuration space $\M^*_{s_i}$. This construction yields a point 

$$Q(X, D) := (Q_2, \dots, Q_n) \in \M = \M^{\tau_2}_{s_2} \times \cdots \times \M^{\tau_n}_{s_n},$$

\noindent where $\tau_i \in \{ +, * \}$ represents the type of the corresponding component $C_i$. $Q(X, D)$ is called the \textit{configuration invariant of} $(X, D)$.
 
Further, performing elementary transformations in $(X, D)$ with centers in $C_0$ does neither change $\Phi_0$ nor the extended divisor (except for the weight $C_1^2$). Hence, it leaves the $s_i$ and $Q(X, D)$ invariant and we can define the configuration invariant for every $m$-standard completion of $V$.\\

\begin{prop}\label{MatchingPrinciple} (Matching Principle, \cite{FKZ4}, Prop. 3.3.1) Let $V = X \backslash D$ be a smooth Gizatullin surface completed by a standard zigzag $D$. Consider the reversed completion $(X^\vee, D^\vee)$ with boundary zigzag $D^\vee = C_0^\vee \cup \cdots \cup C_n^\vee$, the associated numbers $s'_2, \dots, s'_n$ and the types $\tau'_2, \dots, \tau'_n$. Then $s_i = s'_{i^\vee}$ and $\tau_i = \tau'_{i^\vee}$ for all $i = 2, \dots, n$. Moreover, the associated points $Q(X, D)$ and $Q(X^\vee, D^\vee)$ in $\M$ coincide under the natural identification

$$\M = \M^{\tau_2}_{s_2} \times \cdots \times \M^{\tau_n}_{s_n} \cong \M^{\tau'_n}_{s'_n} \times \cdots \times \M^{\tau'_2}_{s'_2}.$$

\end{prop}

\section{Smooth Gizatullin Surfaces with a $(-1)$-completion and Transitivity}

\subsection{Presentations of smooth Gizatullin surfaces}\label{PresentationGizatullinSurface}

\noindent Every smooth Gizatullin surface $V$ can be constructed along with a standard completion $(X, D)$ via a sequence of blow-ups starting from the quadric $Q = \p^1 \times \p^1$. \cite{FKZ4}, 4.1.1, describes an algorithm for constructing standard completions if every boundary component is of type $+$. In essence, it constructs intermediate surfaces along with certain coordinate systems $(x_i, y_i)$, which give the correspondence fibration for a certain pair $(C^\vee_{s^\vee}, C_s)$. In this section we give a generalization of this algorithm. 

On the quadric $Q = \p^1 \times \p^1$ we fix homogeneous coordinates $((s_0: t_0), (s_1: t_1))$ and introduce the affine coordinates $(x_0, y_0)$ via $x_0 := t_0/s_0$ and $y_0 := t_1/s_1$. Furthermore, we let $C_0 = \p^1 \times \{ \infty \}$, $C_1 = \{ \infty \} \times \p^1$ and $C_2 = \p^1 \times \{ 0 \}$. Assume, that among $C_2, \dots, C_n$ the components $C_{k_2} = C_2, C_{k_3}, \dots, C_{k_{r - 1}}, C_{k_r} = C_n$, $k_2 = 2 < k_3 < \dots < k_{r - 1} < k_r = n$, are of type $+$, and that all other components are of type $*$. We choose points $c_{k_3}, \dots, c_{k_r} = c_n$ with $c_{k_i} \in C_{k_{i - 1}} \backslash C_{k_{i - 1} - 1} \cong \a^1$ and finite subsets 

$$M_i \subseteq \begin{cases} C_i \backslash C_{i - 1} \cong \a^1 &, \ i \in \{ k_2, k_3, \dots, k_r \} \\ C_i \backslash (C_{i - 1} \cup C_{i + 1}) \cong \c^* &, \ i \in \{ 2, \dots, n \} \backslash \{ k_2, k_3, \dots, k_r \}. \end{cases}$$

Furthermore, we let $M_i = \{ P_{i, 1}, \dots, P_{i, r_i} \}$. In the coordinates introduced below the subsets $M_i$ are base points of the feathers with mother component $C_i$ and the points $c_{k_i}$ represent coordinates where a $+$-component $C_{k_i}$ is born by an outer blow-up on $C_{k_{i - 1}} \backslash C_{k_{i - 1} - 1}$. Moreover, for a finite subset $M \subseteq \a^1$ or $M \subseteq \c^*$ we abbreviate $P_M(u) := \prod_{a \in M} (u - a) \in \c[u]$.

We consider a decomposition of $X \to Q$ into blow-ups, such that every blow-up either creates a boundary component or a family of feathers attached to the same component:

$$X = X_N \to X_{N - 1} \to \cdots \to X_1 \to X_0 = Q.$$

The following algorithm yields affine coordinate systems $(x_i, y_i)$, $i = 2, \dots, n$, on $X$ such that $C_i = \overline{\{ y_i = 0 \}}$ and $C_{i + 1} = \overline{\{ x_i = 0 \}}$ and such that the correspondence fibration for the pair $(C^\vee_{i^\vee}, C_i)$ is given, in suitable coordinates on $\a^1 = \p^1 \backslash \{ \infty \}$, by the map $(x_i, y_i) \mapsto x_i$ (see Prop. \ref{CorrespondenceFibration} below).

\vspace{15pt}
\noindent \textbf{The algorithm:}

\begin{itemize}
\item[(1)] Let $X_1 \to X_0 = Q$ create the feathers attached to $C_2$, that is, $X_1 \to X_0$ is the blow-up of $M_2 \subseteq C_2 \backslash C_1$. If $P_{2, j}$ has $(x_0, y_0)$-coordinates $(a_{2, j}, 0)$, we let $F_{2, j}$ denote the exceptional curves of the blow-up in $M_2$ and we let $F^\vee_{2, j}$ denote the proper transform on $X_1$ of the closure of the affine line $\{ x_0 = a_{2, j} \}$. Introduce the coordinates 

\begin{equation}\label{Transformation0}
(x_2, y_2) := \left(x_0 - c_{k_3}, \frac{y_0}{P_{M_2}(x_0)} \right).
\end{equation}

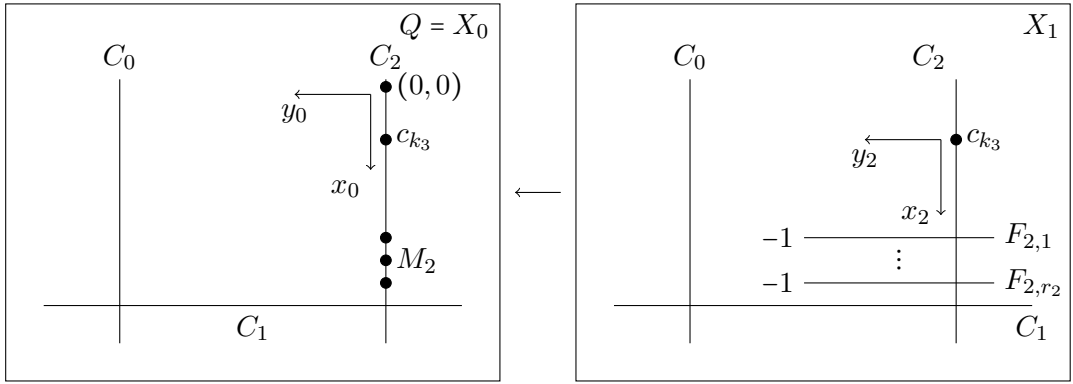
\bigskip
\begin{tikzpicture}
	\draw (0, 0) rectangle (6.5, 5) node[anchor=north east]{$Q = X_0$};
	\begin{scope}
		\draw (1.5, 4) node[anchor=south] {$C_0$} -- (1.5, 0.5);					
		\draw (0.5, 1) -- node[anchor=north] {$C_1$} (6, 1);					
		\draw (5, 4) node[anchor=south] {$C_2$} -- (5, 0.5);					
	\end{scope}	

	\draw [->] (4.8, 3.8) -- (4.8, 2.8) node[anchor=north east] {$x_0$};		\draw [->] (4.8, 3.8) -- (3.8, 3.8) node[anchor=north] {$y_0$};  		
	\filldraw [black] (5, 3.9) circle (2pt) node[anchor=west] {$(0, 0)$};

	\filldraw [black] (5, 3.2) circle (2pt) node[anchor=west] {$c_{k_3}$};
	\filldraw [black] (5, 1.9) circle (2pt);
	\filldraw [black] (5, 1.6) circle (2pt) node[anchor=west] {$M_2$};
	\filldraw [black] (5, 1.3) circle (2pt);

	\draw [->] (7.3, 2.5) -- (6.7, 2.5);

	\draw (7.5, 0) rectangle (14, 5) node[anchor=north east]{$X_1$};
	\begin{scope}
		\draw (9, 4) node[anchor=south] {$C_0$} -- (9, 0.5);					
		\draw (8, 1) -- (13.5, 1) node[anchor=north] {$C_1$};					
		\draw (12.5, 4) node[anchor=south east] {$C_2$} -- (12.5, 0.5);					
	\end{scope}	

	\draw [->] (12.3, 3.2) -- (12.3, 2.2) node[anchor=east] {$x_2$};			
	\draw [->] (12.3, 3.2) -- (11.3, 3.2) node[anchor=north] {$y_2$};  	

	\filldraw [black] (12.5, 3.2) circle (2pt) node[anchor=west] {$c_{k_3}$};
	\draw (13, 1.9) node[anchor=west] {$F_{2, 1}$} -- node[anchor=north] {$\vdots$} (10.5, 1.9)     node[anchor=east] {$-1$};					
	
	\draw (13, 1.3) node[anchor=west] {$F_{2, r_2}$} -- (10.5, 1.3) node[anchor=east] {$-1$};					
\end{tikzpicture}

\bigskip 
\item[(2)] Let the blow-up $X_2 \to X_1$ in $c_{k_3} \in C_2 \backslash C_1 \cong \a^1$ create the component $C_{k_3}$, the first $+$-component after $C_2$. Introduce the coordinates $(w_2, z_2)$ and $(w_{k_3}, z_{k_3})$ by 

$$(w_2, z_2) := \left( x_2, \frac{y_2}{x_2} \right) \quad \text{and} \quad (w_{k_3}, z_{k_3}) := \left( \frac{x_2}{y_2}, y_2 \right).$$

\noindent The situation is illustrated in the figure below (the red dashed line is the proper transform of $\{ x_2 = 0 \} = \{ x_0 = c_{k_3} \}$).

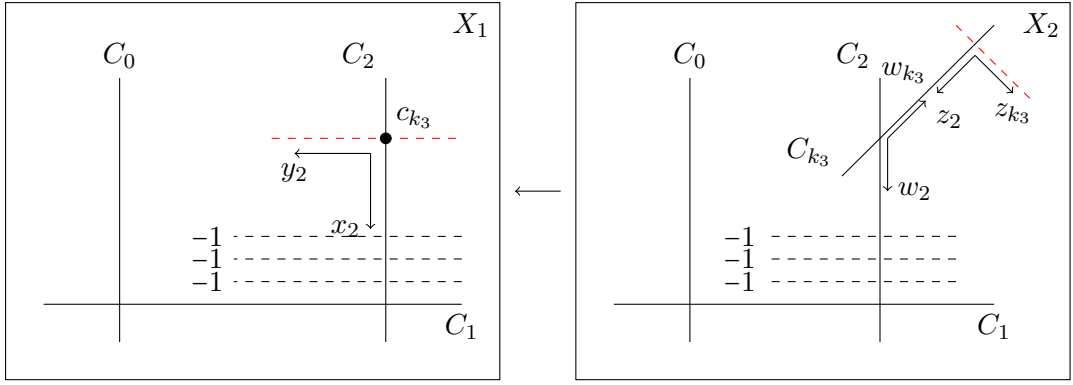
\bigskip
\begin{tikzpicture}
	\draw (0, 0) rectangle (6.5, 5) node[anchor=north east]{$X_1$};
	\begin{scope}
		\draw (1.5, 4) node[anchor=south] {$C_0$} -- (1.5, 0.5);					
		\draw (0.5, 1) -- (6, 1) node[anchor=north] {$C_1$};					
		\draw (5, 4) node[anchor=south east] {$C_2$} -- (5, 0.5);					
	\end{scope}	

	\draw [->] (4.8, 3) -- (4.8, 2) node[anchor=east] {$x_2$};			
	\draw [->] (4.8, 3) -- (3.8, 3) node[anchor=north] {$y_2$};  		

	\draw[dashed, color=red] (3.5, 3.2) -- (6, 3.2);					

	\filldraw [black] (5, 3.2) circle (2pt) node[anchor=south west] {$c_{k_3}$};
	\draw[dashed] (6, 1.9) -- (3, 1.9) node[anchor=east] {$-1$};					
	\draw[dashed] (6, 1.6) -- (3, 1.6) node[anchor=east] {$-1$};					
	\draw[dashed] (6, 1.3) -- (3, 1.3) node[anchor=east] {$-1$};

	\draw [->] (7.3, 2.5) -- (6.7, 2.5);

	\draw (7.5, 0) rectangle (14, 5) node[anchor=north east]{$X_2$};
	\begin{scope}
		\draw (9, 4) node[anchor=south] {$C_0$} -- (9, 0.5);					
		\draw (8, 1) -- (13, 1) node[anchor=north] {$C_1$};					
		\draw (11.5, 4) node[anchor=south east] {$C_2$} -- (11.5, 0.5);					
		\draw (11, 2.7) node[anchor=south east] {$C_{k_3}$} -- (13, 4.7);					
	\end{scope}	

	\draw[dashed, color=red] (12.5, 4.7) -- (13.5, 3.7);					

	\draw [->] (11.6, 3.2) -- (11.6, 2.5) node[anchor=west] {$w_2$};			
	\draw [->] (11.6, 3.2) -- (12.1, 3.7) node[anchor=north west] {$z_2$};  		

	\draw [->] (12.75, 4.3) -- (12.25, 3.8) node[anchor=south east] {$w_{k_3}$};			
	\draw [->] (12.75, 4.3) -- (13.25, 3.8) node[anchor=north] {$z_{k_3}$};  		

	\draw[dashed] (12.5, 1.9) -- (10, 1.9) node[anchor=east] {$-1$};					
	\draw[dashed] (12.5, 1.6) -- (10, 1.6) node[anchor=east] {$-1$};					
	\draw[dashed] (12.5, 1.3) -- (10, 1.3) node[anchor=east] {$-1$};					
\end{tikzpicture}

\bigskip
\noindent Now we perform inner blow-ups to create the $*$-components $C_3, \dots, C_{k_3 - 1}$ and we fix some $t \in \{ 2, \dots, k_3 - 1 \}$. We start with the blow-up in $C_2 \cap C_{k_3}$ and we assume by induction that some inner components as well as the corresponding affine coordinates $(w_i, z_i)$ are already created, and that $C' = \overline{\{ z_i = 0 \}}$ and $C'' = \overline{\{ w_i = 0 \}}$, where $C'$ and $C''$ are two irreducible neighboring components of the zigzag such that $C'$ precedes $C''$. Blowing up the intersection point $C' \cap C''$, we introduce new affine coordinates $(w_{i + 1}, z_{i + 1})$ either via 

\begin{equation}\label{TwoTypesOfTransformations}
(w_{i + 1}, z_{i + 1}) = \left( \frac{w_i}{z_i}, z_i \right) \quad \text{or via} \quad (w_{i + 1}, z_{i + 1}) = \left( w_i, \frac{z_i}{w_i} \right).
\end{equation}

\noindent Denoting by $E$ the exceptional curve of the last blow-up, we obtain $E = \overline{\{ z_{i + 1} = 0 \}}, C'' = \overline{\{ w_{i + 1} = 0 \}}$ in the first case and $C' = \overline{\{ z_{i + 1} = 0 \}}, E = \overline{\{ w_{i + 1} = 0 \}}$ in the second case. In the end, after choosing an appropriate ordering of transformations as in (\ref{TwoTypesOfTransformations}), after the last blow-up $X_{k_3 - 2} \to X_{k_3 - 3}$ we have $C_t = \{ z_{k_3} = 0 \}$ and $C_{t + 1} = \{ w_{k_3} = 0 \}$. We denote these resulting coordinates by $(w_t, z_t)$ instead of $(w_{k_3}, z_{k_3})$.

\bigskip
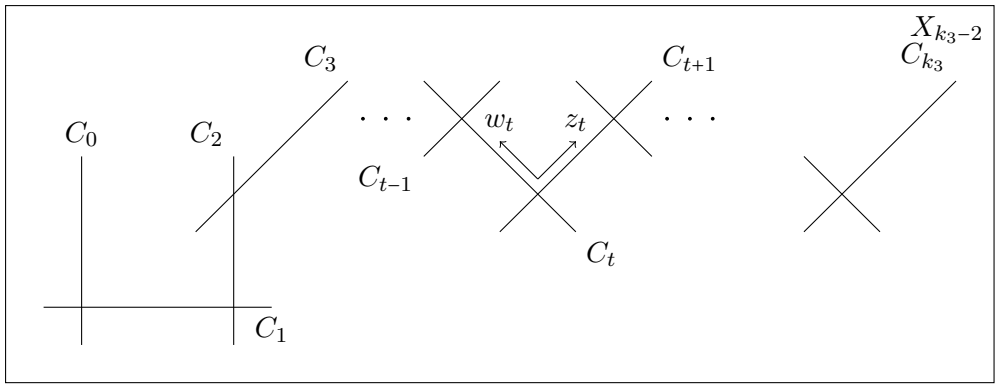
\begin{tikzpicture}
	\draw (0, 0) rectangle (13, 5) node[anchor=north east]{$X_{k_3 - 2}$};
	\begin{scope}
		\draw (1, 3) node[anchor=south] {$C_0$} -- (1, 0.5);					
		\draw (0.5, 1) -- (3.5, 1) node[anchor=north] {$C_1$};					
		\draw (3, 3) node[anchor=south east] {$C_2$} -- (3, 0.5);					
 		\draw (2.5, 2) -- (4.5, 4) node[anchor=south east] {$C_3$};					
    	\filldraw [black] (4.7, 3.5) circle (0.5pt);
    	\filldraw [black] (5, 3.5) circle (0.5pt);
    	\filldraw [black] (5.3, 3.5) circle (0.5pt);
		\draw (5.5, 3) node[anchor=north east] {$C_{t - 1}$} -- (6.5, 4);					
		\draw (5.5, 4) -- (7.5, 2) node[anchor=north west] {$C_t$};					
		\draw (6.5, 2) -- (8.5, 4) node[anchor=south west] {$C_{t + 1}$};					
		\draw (7.5, 4) -- (8.5, 3);					
    	\filldraw [black] (8.7, 3.5) circle (0.5pt);
    	\filldraw [black] (9, 3.5) circle (0.5pt);
    	\filldraw [black] (9.3, 3.5) circle (0.5pt);
		\draw (10.5, 3) -- (11.5, 2);					
		\draw (10.5, 2) -- (12.5, 4) node[anchor=south east] {$C_{k_3}$};					
	\end{scope}	

	\draw [->] (7, 2.7) -- (6.5, 3.2) node[anchor=south] {$w_t$};			
	\draw [->] (7, 2.7) -- (7.5, 3.2) node[anchor=south] {$z_t$};  		
\end{tikzpicture}

\bigskip
\item[(3)] Performing coordinate transformations as in (2), it is easy to check that we obtain some coordinate systems $(w_i, z_i)$, $2 \leq i \leq k_3$, with $C_i = \overline{\{ z_i = 0 \}}$ and $C_{i + 1} = \overline{\{ w_i = 0 \}}$ (if $i + 1 \leq k_3$), which are related by 

\begin{equation}\label{Transformation1}
(w_j, z_j) = T_{ij}(w_i, z_i) := (w_i^{k_{ij}}z_i^{l_{ij}}, w_i^{p_{ij}}z_i^{q_{ij}}), \quad \left|\begin{array}{cc} k_{ij} & p_{ij} \\ l_{ij} & q_{ij} \end{array}\right| = 1, \ q_{ij} > 0, \ l_{ij} < 0 \quad \forall \ j > i. 
\end{equation}

\noindent Moreover, $k_{ij} = 0$, $l_{ij} = -1$ holds if and only if $j = i + 1$ (\ie $w_{i + 1} = z_i^{-1}$).
\item[(4)] For each $i = 3, \dots, k_3 - 1$, we let $P_{i, j}$ have the $(w_i, z_i)$-coordinates $(a_{i, j}, 0)$. Let $F_{i, j}$, $j \in \{ 1, \dots, r_i \}$, denote the exceptional curves of the blow-up in $M_i$ and let $F^\vee_{i, j}$ denote the proper transform of the closure of the affine line $\{ w_i = a_{i, j} \}$. Introduce the coordinates

\begin{equation}\label{Transformation2}
(x_i, y_i) := T^{\text{Bl.up}}_i(w_i, z_i) := \left( w_i, \frac{z_i}{P_{M_i}(w_i)} \right).
\end{equation}

\noindent Moreover, replace the coordinates $(w_j, z_j)$ by $(T_{ij} \circ T^{\text{Bl.up}}_i \circ T^{-1}_{ij})(w_j, z_j)$ for all $j = i + 1, \dots, k_3$. In particular, using (3) it follows that $w_{i + 1} = y_i^{-1}$.
\item[(5)] Assume now by induction that all components $C_2, \dots, C_{k_s}$, $2 \leq s < r$, including the feathers with mother component $C_i$, $i \leq k_s - 1$, are already created. We repeat the steps (1) - (4), but now with the coordinates $(w_{k_s}, z_{k_s})$ instead of $(x_0, y_0)$.
\item[(6)] Finally, if all feathers on $C_2, \dots, C_{n - 1}$ are created, we create in the same way feathers on $C_n$ and introduce similarly the coordinates $(x_n, y_n)$ and the curves $F_{n, j}$ and $F^\vee_{n, j}$, $j = 1, \dots, r_n$. 
\end{itemize}

\noindent We denote the smooth projective surface $X$ obtained in this way by 

$$X = X_n := X(M_2, M_3, \dots, M_{k_3 - 1}, c_{k_3}, M_{k_3}, \dots, M_{k_4 - 1}, c_{k_4}, M_{k_4}, \dots, c_n, M_n)$$

\noindent and the boundary divisor is

$$D := C_0 \cup \cdots \cup C_n.$$

By construction, the surface $X \backslash D = V$ is a smooth Gizatullin surface. Note, that $X \backslash D$ is affine: by applying appropriate elementary transformations we may assume that $C^2_1 \gg 0$ and hence it is the support of an ample divisor by the Nakai-Moishezon criterion. We refer to $X_n$ as a \emph{presentation} of $V$.

\begin{df}\label{DfMinusOneCompletion} A presentation $X_n$ of a smooth Gizatullin surface $V$ is called a $(-1)$-presentation (or of $(-1)$-type) if $c_{k_i} \not\in M_{k_{i - 1}}$ for all $i = 3, \dots, r$. Equivalently, all feathers in $X_n$ have self-intersection $-1$.
 
Similarly, a standard completion $(X, D)$ of $V$ is called a $(-1)$-completion if all feathers are $(-1)$-feathers, or equivalently, if all feathers are attached to their mother components.

Obviously, $X_n$ gives a $(-1)$-presentation of $V$ if and only if $(X_n, D)$ is a $(-1)$-completion of $V$.
\end{df}

The coordinate systems $(x_i, y_i)$ introduced in the above algorithm reflects the correspondence fibration. Two particular cases of the following proposition were shown in \cite{FKZ4}, Prop. 5.2.1 (if all components of $D$ are outer components) and in \cite{Ko}, Prop. 2.36 (if the components $C_3, \dots, C_{n - 1}$ are inner components):

\begin{prop}\label{CorrespondenceFibration} Assume that there are given affine coordinates $(x_i, y_i)$ on a presentation $X_n$ as in \ref{PresentationGizatullinSurface}. Then, in appropriate coordinates on $\a^1 = \p^1 \backslash \{ \infty \}$, the map $\overline{\pi}_{j - 1} := \Phi_{|C^\vee_{j^\vee - 1}|}: W \to \p^1$, being the correspondence fibration for the pair $(C^\vee_{j^\vee}, C_j)$, is given by $x_j$. In particular, the pair $(F_{i, j}, F^\vee_{i, j})$ is a pair of matching feathers.
\end{prop}

\begin{bew} We show this by induction on the number $r$ of outer components of $D$. For $r = 2$ the claim is precisely the one of Porp. 2.36 in \cite{Ko}. Assume now that the claim holds for all $j = 2, \dots, k_r = n$. We create one new outer component $C_{k_{r + 1}}$ with coordinates $(w_{k_{r + 1}}, z_{k_{r + 1}})$. If $k_{r + 1} = k_r + 1$, that is, there are no $*$-components between $C_{k_r}$ and $C_{k_{r + 1}}$, then we create all feathers on $C_{k_{r + 1}}$. In this case the claim follows by \cite{FKZ4}, Prop. 5.2.1, since the inductive construction of the coordinates in the algorithm \cite{FKZ4}, 5.1.1 do not differ from the one in \ref{PresentationGizatullinSurface}.

Assume now, that $k_{r + 1} > k_r + 1$. We create inner components $C_{k_r + 1}, \dots, C_{k_{r + 1} - 1}$ together with all feathers which are attached to these components. Knowing by induction hypothesis that $(x_{k_r}, y_{k_r}) \mapsto x_{k_r}$ gives the correspondence fibration for the pair $(C^\vee_{{k_r}^\vee}, C_{k_r})$, it follows again by \cite{Ko}, Prop. 2.36 that $(x_j, y_j) \mapsto x_j$ gives the correspondence fibration for the pair $(C^\vee_{j^\vee}, C_j)$ for all $j = k_r + 1, \dots, k_{r + 1} - 1$.

Now, at this step there are no feathers attached to $C_{k_{r + 1}}$. The map $(w_{k_{r + 1}}, z_{k_{r + 1}}) \mapsto w_{k_{r + 1}}$ gives the correspondence fibration for the pair $(C^\vee_{k_{r + 1}^\vee}, C_{k_{r + 1}})$ since our coordinates are related by $w_{k_{r + 1}} = \frac{1}{y_{k_{r + 1} - 1}}$ (keep in mind that the coordinates $(w_{k_{r + 1}}, z_{k_{r + 1}})$ change while creating inner components, see step (4) in the above algorithm). Finally, creating feathers on $C_{k_{r + 1}}$ does not neither change the coordinate $w_{k_{r + 1}}$ nor the correspondence firbration for $(C^\vee_{k_{r + 1}^\vee}, C_{k_{r + 1}})$. After creating all feathers the coordinate $w_{k_{r + 1}}$ becomes $x_{k_{r + 1}}$ and the proof is complete.
\end{bew}

The group 

$$\Jon = \Aut_{x_0}(\a^2) = \{ (x_0, y_0) \mapsto (ax_0 + P(y_0), by_0) \mid a, b \in \c^*, P(y_0) \in \c[y_0] \}$$

\noindent acts on presentations $X_n$ in the following way. Given an element $h \in \Jon$, the set $M_2$ is moved by $h$ into a new set of points $M'_2 \subseteq C_2$. Thus, $h$ induces an isomorphism $X(M_2) \backslash C_0 \stackrel{\sim}{\to} X(M'_2) \backslash C_0$. Under this isomorphism $c_{k_3}$ is mapped to a point $c'_{k_3}$, yielding an isomorphism $X(M_2, c_{k_3}) \backslash C_0 \stackrel{\sim}{\to} X(M'_2, c'_{k_3}) \backslash C_0$. Continuing this way (for outer and inner components) we obtain an isomorphism

\begin{eqnarray*}
&& X(M_2, M_3, \dots, M_{k_3 - 1}, c_{k_3}, M_{k_3}, \dots, c_n, M_n) \backslash C_0 \\
&\stackrel{\sim}{\to}& X(M'_2, M'_3, \dots, M'_{k_3 - 1}, c'_{k_3}, M'_{k_3}, \dots, c'_n, M'_n) \backslash C_0.
\end{eqnarray*}

\noindent We denote this new presentation of $V$ by 

$$h.X(M_2, M_3, \dots, M_{k_3 - 1}, c_{k_3}, M_{k_3}, \dots, c_n, M_n).$$

\noindent Note that to every standard completion $(X, D)$ of $V$ we can associate a presentation $X_n$ of $V$. Therefore we can write $h.(X, D)$ or $(h.X, h.D)$ instead of $h.X_n$.

Note also, that $h$ maps the set of points blown up by $(X, D) \to Q$ onto the set of points blown up by $h.(X, D) \to Q$. Thus, $h$ can be considered as a fibered map $h: (X, D) \dasharrow h.(X, D)$. In particular, the image $h(F)$ of a feather $F$ of $D_{\text{ext}}$ is a feather of $(h.D)_{\text{ext}}$.

Now, every $h \in \Jon$ can be decomposed into $h = t \circ \varphi_k \circ \cdots \circ \varphi_1$, such that $t \in \t \cong (\c^*)^2$ and $\varphi_i$ is of the form $\varphi_i(x_0, y_0) = (x_0 + a_iy_0^{t_i - 2}, y_0)$, $a_i \in \c$, $t_i \geq 2$. These $\varphi_i$ are called \emph{elementary shifts}. Thus, for studying the action of $\Jon$ on presentations of $V$ it suffices to study the actions of the elementary shifts and of $\t$.

\bigskip
\subsection{Transitivity of the automorphism group}

Besides the surface $\c^* \times \c^*$, the quasi-homogeneous surfaces are precisely the surfaces which are completable by a zigzag. Unfortunately, the last property is far away from being equivalent to homogeneity of the surface. In each positive characteristic $p > 0$ counterexamples were found very soon (see \cite{DG1}). But in characteristic zero a long time no counterexamples were known. In \cite{Ko} it is proved that a certain subclass of Gizatullin surfaces which admit a distinguished and rigid extended divisor provides the desired counterexamples (for the notion of a distinguished and rigid extended divisor see \cite{FKZ3} section 1 and 2). But it is an open problem (in every chracteristic) to classify \emph{all} quasi-homogeneous surfaces which are not homogeneous.

In this section we provide a partial solution for this problem. To be more precise, we give conditions for a quasi-homogeneous surface to be homogeneous if it admits a $(-1)$-completion.\footnote[1]{Examples for quasi-homogeneous surfaces which do not admit a $(-1)$-completion are not known yet. But the author presumes their existence.}

First, let us recall the "canonical" candidates for those points which might be contained in the finite set $V \backslash O$:

\begin{prop}\label{AutFiniteComplement}(\cite{Ko}, Prop. 3.7) Let $V$ be a smooth Gizatullin surface with a standard completion $(X, D)$ and the associated $\p^1$-fibration $\overline{\pi}: X \to \p^1$. Furthermore, let $F_{i, \rho}$ be the feathers of the extended divisor $D_{\text{ext}}$ which have mother component $C_i$, and let $F^\vee_{i, \rho}$ be the matching feathers in $D^\vee_{\text{ext}}$. Letting $O$ be big orbit of the natural action of $\Aut(V)$, we have

$$V \backslash O \subseteq \bigcup_{i, j} (F_{i, \rho} \cap F^\vee_{j, \sigma}) \footnote[2]{Note, that in the case i = j the intersection $F_{i, \rho} \cap F^\vee_{i, \sigma}$ is non-empty if and only if $\rho = \sigma$.}.$$

\end{prop}

Proposition \ref{AutFiniteComplement} shows that the finite complement $V\backslash O$ of the big orbit $O$ is contained in the union of all feathers. In order to formulate the main result, we need to determine those boundary components $C_i$ with the property that for every feather $F$, which is attached to $C_i$, we have $F \backslash D \subseteq O$.

\begin{df}\label{DfExceptionalComponent} Let $X = X(M_2, \dots, M_{k_3 - 1}, c_{k_3}, M_{k_3}, \dots, c_n, M_n)$ be a presentation of a smooth Gizatullin surface $V$. Then $X$ arises as a blowup of the quadric $Q = \p^1 \times \p^1$ and, for each $s = 3, \dots, r - 1$, a suitable ordering of the blow-ups yields an intermediate surface $(X', D')$, such that the dual graph of $D'$ has the form

$$\dots \lin \cou{C_{k_s}}{w'} \lin \cou{C_{\tau_{s, 1}}}{-2} \lin \dots \lin \cou{C_{\tau_{s, m_s - 1}}}{-2} \lllin \cou{C_{\tau_{s, m_s}}}{-1} \lin \cou{C_{k_{s + 1}}}{w''} \lin \dots$$

\noindent such that $m_s$ is maximal. The proper transform of $C_{\tau_{s, 1}}, \dots, C_{\tau_{s, m_s}}$ under the blowup $X \to X'$ are called \emph{exceptional components of} $D$. We denote $\mathfrak{E}_D := \{ \tau_{s, j} \mid s \in \{ 2, \dots, r - 1 \}, j \in \{ 1, \dots, m_s \} \}$, and further, we denote the set $\{ \tau^\vee_{s, j} \mid s \in \{ 2, \dots, r - 1 \}, j \in \{ 1, \dots, m_s \} \}$ by $\mathfrak{E}^\vee_D$. Note, that these $C_{\tau_{s, i}}$ are uniquely determined by this condition.
\end{df}

For example, if $D_{\text{ext}}$ has the dual graph

\vspace{15pt}
$$D_{\text{ext}}: \quad \cu{0} \lin \cu{0} \lin \cu{-3} \lin \cu{-2} \nlin \cshiftup{-1}{} \lin \cu{-3} \lin \cu{-2} \nlin \cshiftup{-1}{} \lin \cu{-3}$$

\vspace{15pt}
\noindent then $C_4$ and $C_5$ are the exceptional components of $D$. Indeed, the blowup $X \to Q$ factorizes as follows:

\vspace{15pt}
$$\cou{C_0}{0} \lin \cou{C_1}{0} \lin \cou{C_2}{0} \quad \leftarrow \quad \cou{C_0}{0} \lin \cou{C_1}{0} \lin \cou{C_2}{-2} \lin \cou{C_4}{-2} \lin \cou{C_5}{-1} \lin \cou{C_6}{-3} \quad \leftarrow \quad \cu{0} \lin \cu{0} \lin \cu{-3} \lin \cu{-2} \nlin \cshiftup{-1}{} \lin \cu{-3} \lin \cu{-2} \nlin \cshiftup{-1}{} \lin \cu{-3}$$

\vspace{15pt}
\noindent Hence $\mathfrak{E}_D = \{ 4, 5 \}$.

\vspace{15pt}
Most of the technical work is contained in the following lemma:

\begin{lemma}\label{IntersectionPointOfFeathersByShifts} Let $V$ be a smooth Gizatullin surface which admits a $(-1)$-completion and let $(X, D)$ be a standard completion of $V$, corresponding to a presentation $X_n = X(M_2, \dots, M_{k_3 - 1}, c_{k_3}, M_{k_3}, \dots, c_n, M_n)$. Consider an elementary shift $h(x_0, y_0) = (x_0 + ay^{t - 2}_0, y_0)$, $a \in \c$ and $t \geq 2$, and the corresponding map $h: (X, D) \dasharrow (X', D') := h.(X, D)$. Given a feather $F$ of $D_{\text{ext}}$, we denote by $F^\vee$ its matching feather in $D^\vee_{\text{ext}}$. Let $C_i$ be an inner component of $D$, which is not an exceptional one. If $F$ is a feather of $D_{\text{ext}}$ attached to $C_i$ and $G$ its image under $h$, then  

$$h(F \cap F^\vee) = G \cap G^\vee.$$

\end{lemma}

\begin{bew} We consider on $(X, D)$ the coordinates $(x_i, y_i)$, $2 \leq i \leq n$, which were constructed in \ref{PresentationGizatullinSurface} and define 

$$(\xi_i, \eta_i) := h(x_i, y_i).$$

\noindent Sometimes we will consider the lift of $h$ on an intermediate surface with our coordinates $(w_i, z_i)$. In this case we also write 

$$(\xi_i, \eta_i) = h(w_i, z_i),$$

\noindent if no confusion arises. Now we are interested in the behaviour of $h$ in the infinitely near neighborhood of $C_2$. Since $\xi_i$ and $\eta_i$ are rational functions in $x_i, y_i$ which are regular in a neighborhood of $C_2$, we can express them as power series in $x_i, y_i$. This will cause less difficulties to determine the general \emph{form} of the lift of $h$ after performing several blow-ups.\\

\noindent \textbf{Claim 1:} The elementary shift $h$ is a rational function in $x_i, y_i$ (respectively $w_i, z_i$), expressed by power series in $\c[[x_i, y_i]]$ (respectively $\c[[w_i, z_i]]$), and is given by the following expressions:

\begin{itemize}
\item[(1)] If $i = k_s$ and $2 \leq  s \leq t - 1$, then

\begin{eqnarray*}
\xi_i &=& w_i + z_i^{t - s}R'(w_i, z_i),\\
\eta_i &=& z_i(1 + z_i^{t - s}S'(w_i, z_i)).
\end{eqnarray*}

\noindent for some power series $R'(w_i, z_i), S'(w_i, z_i) \in \c[[w_i, z_i]]$. In particular, the same holds for the coordinates $(x_i, y_i)$, \ie

\begin{eqnarray*}
\xi_i &=& x_i + y_i^{t - s}R(x_i, y_i),\\
\eta_i &=& y_i(1 + y_i^{t - s}S(x_i, y_i)),
\end{eqnarray*}

\noindent for some power series $R$ and $S$. 
\item[(2)] If $i = k_s$ and $t \leq  s \leq r$, then

\begin{eqnarray*}
\xi_i &=& w_i + f_i(a) + z_iR'(w_i, z_i),\\
\eta_i &=& z_i(1 + z_iS'(w_i, z_i)),
\end{eqnarray*}

\noindent for some power series $R'(w_i, z_i), S'(w_i, z_i) \in \c[[w_i, z_i]]$. In particular, the same holds for the coordinates $(x_i, y_i)$, \ie

\begin{eqnarray*}
\xi_i &=& x_i + f_i(a) + y_iR(x_i, y_i),\\
\eta_i &=& y_i(1 + y_iS(x_i, y_i)),
\end{eqnarray*}

\noindent for some constant $f_i(a) \in \c$, depending on $a$, and some $R(x_i, y_i), S(x_i, y_i) \in \c[[x_i, y_i]]$.
\item[(3)] Consider the lift $(\xi_i, \eta_i)$ of $h$ in the coordinates $(w_i, z_i)$ defined in \ref{PresentationGizatullinSurface} step (2) - (4). If $k_s < i < k_{s + 1}$ for some $s = 2, \dots, r - 1$, then

\begin{eqnarray*}
\xi_i &=& w_i(1 + w_i^{a_i}z_i^{b_i}R(w_i, z_i)),\\
\eta_i &=& z_i(1 + w_i^{a_i}z_i^{b_i}S(w_i, z_i)),
\end{eqnarray*}

\noindent for some integers $a_i \geq 0$, $b_i \geq 1$.  Moreover, we have $b_i \geq 2$ if and only if $C_i$ is not an exceptional component.
\end{itemize}

We show Claim 1 by induction on the number of blow-ups performed to obtain $X_n$ from $Q$. For this, we have to check that the general form of $(\xi_i, \eta_i)$ is preserved under\\

\noindent (i) inner blow-ups, \\
\noindent (ii) outer blow-ups, \\
\noindent (iii) creating feathers on outer boundary components, \\
\noindent (iv) creating feathers on inner boundary components. \\

First, the claim obviously holds for the surface $Q$. Assume now that the intermediate surface $X(M_2, \dots, c_{k_s}, M_{k_s})$, $s \leq t - 1$ is already constructed. We use algorithm \ref{PresentationGizatullinSurface} and create an outer component $C_{k_{s + 1}}$ by a blow-up in a point in $C_{k_s} \backslash C_{k_s - 1}$, and thereafter, we create the inner components $C_{k_s + 1}, \dots, C_{k_{s + 1} - 1}$. Now, if $k_s < i < k_{s + 1} - 1$, the lift of $h$ has the form

\begin{equation}\label{ShiftInner}
\xi_i = w_i(1 + w_i^{a_i}z_i^{b_i}R(w_i, z_i)), \quad \eta_i = z_i(1 + w_i^{a_i}z_i^{b_i}S(w_i, z_i)),
\end{equation}

\noindent with some integers $a_i \geq 0$, $b_i \geq 1$ some power series $R$ and $S$. Moreover, we have $b_i \geq 2$ if and only if $C_i$ is not an exceptional component. This can be shown in the same way as in the proof of \cite{Ko}, Lemma 3.12. At this stage, we did not create any feather on the inner components $C_{k_s + 1}, \dots, C_{k_{s + 1} - 1}$.

Now we start to create feathers in the inner components $C_{k_s + 1}, \dots, C_{k_{s + 1} - 1}$ and show that replacing $(w_i, z_i)$ by $(T_{ji} \circ T^{\text{Bl.up}}_j \circ T^{-1}_{ji})(w_i, z_i)$ and $(\xi_i, \eta_i)$ by $(T_{ji} \circ T^{\text{Bl.up}}_j \circ T^{-1}_{ji})(\xi_i, \eta_i)$ for some $j \in \{ k_s + 1, \dots, i - 1 \}$ preserves the form of equation (\ref{ShiftInner}) and even the values for $a_i$ and $b_i$. The situation can be reduced to the case where $T^{\text{Bl.up}}_j$ creates a single feather at a point $c$ on $(C_j \backslash C_{j - 1} \cup C_{j + 1}) \cong \c^*$. Let as in \ref{PresentationGizatullinSurface} (3) be 

\begin{equation}\label{Transition}
(\xi_i, \eta_i) = T_{ji}(\xi_j, \eta_j) = (\xi_j^k\eta_j^l, \xi_j^p\eta_j^q) \quad \text{with} \quad kq - lp = 1, \ -l, q > 0. 
\end{equation}

\noindent For brevity we let $m := -l > 0$, $a := a_i$, $b := b_i$, $(w', z') := (T_{ji} \circ T^{\text{Bl.up}}_j \circ T^{-1}_{ji})(w_i, z_i)$ and $(\xi', \eta') := (T_{ji} \circ T^{\text{Bl.up}}_j \circ T^{-1}_{ji})(\xi_i, \eta_i)$. Since $T^{\text{Bl.up}}_j(w_j, z_j) = \left( w_j, \frac{z_j}{w_j - c} \right)$ and $T^{\text{Bl.up}}_j(\xi_j, \eta_j) = \left( \xi_j, \frac{\eta_j}{\xi_j - c} \right)$ (note, that $h$ induces the identity on every inner component), we see that

\begin{eqnarray*}
(\xi', \eta') &=& (T_{ji} \circ T^{\text{Bl.up}}_i \circ T^{-1}_{ji})(\xi_i, \eta_i) = \left( \xi_i(\xi_i^q\eta_i^m - c)^m, \frac{\eta_i}{(\xi_i^q\eta_i^m - c)^q} \right),\\
(w', z') &=& \left( w_i(w_i^qz_i^m - c)^m, \frac{z_i}{(w_i^qz_i^m - c)^q} \right).
\end{eqnarray*}

\noindent The inverse relation is $(w_i, z_i) = \left( \frac{w'}{(w'^qz'^m - c)^m}, z'(w'^qz'^m - c)^q \right)$. Putting this together and denoting by $R, S, T$ etc. the resulting power series in the intermediate steps, we get 

\begin{eqnarray*}
\xi' &=& \xi_i(\xi^q_i\eta^m_i - c)^m \\
&=& w_i(1 + w_i^az_i^bR)(w_i^qz_i^m(\underbrace{1 + w_i^az_i^bR)(1 + w_i^az_i^bS)}_{=: 1 + w_i^az_i^bT} - c)^m \\
&=& \frac{w'}{({w'}^q{z'}^m - c)^m}(1 + {w'}^a{z'}^b\tilde{R}) \cdot [({w'}^q{z'}^m - c) + {w'}^{q + a}{z'}^{q + b}\underbrace{({w'}^q{z'}^m - c)^{bq - as}T}_{=: U}]^m
\end{eqnarray*}

\noindent Note, that $(w^qz^m - c)^{bq - as}$ admits a power series expansion near $z' = 0$ and therefore $U := (w^qz^m - c)^{bq - as}T$ is a power series as well. Using binomial expansion we obtain

\begin{eqnarray*}
\xi' &=& \frac{w'}{({w'}^q{z'}^m - c)^m}(1 + {w'}^a{z'}^b\tilde{R}) \cdot \left[ ({w'}^q{z'}^m - c)^m + \sum_{\mu = 1}^m \binom{m}{\mu} ({w'}^q{z'}^m - c)^{m - j} {w'}^{\mu(q + a)}{z'}^{\mu(q + b)}U^\mu \right] \\
&=& w'(1' + {w'}^a{z'}^b\tilde{R}) \cdot \left[ 1 + \sum_{\mu = 1}^m \binom{s}{\mu} ({w'}^q{z'}^m - c)^{- j} {w'}^{\mu(q + a)}{z'}^{\mu(q + b)}U^\mu \right]. \\
\end{eqnarray*}

\noindent Now, every exponent of $w'$ and $z'$, respectively, in the expression inside the bracket is at least $a + 1$ and $b + 1$, respectively. Thus, since $({w'}^q{z'}^m - c)^{- j}$ also admits a power series representation, the bracket is of the form $1 + w'^az'^bV$ for some power series $V$. In summary, we get 

$$\xi' = w'(1 + {w'}^a{z'}^bR')$$

\noindent for some power series $R'$. Similarly, $\eta'$ can be computed in the same way. This gives assertion (3) in the case $s \leq t - 1$.

Let us show assertion (1). Assume by induction that (1) holds for some $i = k_s$ such that $s \leq t - 1$. Applying our algorithm in \ref{PresentationGizatullinSurface}, we create $C_{k_{s + 1}}$ by an outer blow-up in a point of $C_{k_s} \backslash C_{k_s - 1}$, say in $(x_{k_s}, y_{k_s}) = (c, 0)$. Since $s \leq t - 1$, we have $h(c, 0) = (c, 0)$. This leads to 

$$(x_{k_s}, y_{k_s}) = \left( w_{k_{s + 1}}z_{k_{s + 1}} + c, z_{k_{s + 1}} \right) \quad \text{and} \quad (\xi_{k_s}, \eta_{k_s}) = \left( \xi_{k_{s + 1}}\eta_{k_{s + 1}} + c, \eta_{k_{s + 1}} \right).$$

\noindent Now, using that $(\xi_{k_s}, \eta_{k_s}) = (x_{k_s} + y_{k_s}^{t - s}R, y_{k_s}(1 + y_{k_s}^{t - s}S))$, it is easy to check that $(\xi_{k_{s + 1}}, \eta_{k_{s + 1}})$ is of the form

$$(\xi_{k_{s + 1}}, \eta_{k_{s + 1}}) = (x_{k_{s + 1}} + y_{k_{s + 1}}^{t - s - 1}R', y_{k_{s + 1}}(1 + y_{k_{s + 1}}^{t - s - 1}S'))$$

\noindent with some power series $R', S'$, if $s \leq t - 2$. Moreover, in the case $s = t - 1$, we even obtain that $(\xi_{k_{s + 1}}, \eta_{k_{s + 1}})$ is of the form

$$(\xi_{k_{s + 1}}, \eta_{k_{s + 1}}) = (x_{k_{s + 1}} + f(a) + R', y_{k_{s + 1}}(1 + y_{k_{s + 1}}S'))$$

\noindent with some constant $f(a) \in \c$, depending on $a$, and some power series $R', S'$ (this arises due to the fact that for any power series $R(x, y)$, the expression $R(uv + c, v)$ is always of the form $b + vS(u, v)$ for some constant $b$ and some power series $S$).

As before, we start to create the inner components $C_{k_s + 1}, \dots, C_{k_{s + 1} - 1}$, lying between the outer components $C_{k_s}$ and $C_{k_{s + 1}}$. Let us now check that replacing $(w_{k_{s + 1}}, z_{k_{s + 1}})$ by $(T_{j, k_{s + 1}} \circ T^{\text{Bl.up}}_j \circ T^{-1}_{j, k_{s + 1}})(w_{k_{s + 1}}, z_{k_{s + 1}})$ and $(\xi_{k_{s + 1}}, \eta_{k_{s + 1}})$ by $(T_{j, k_{s + 1}} \circ T^{\text{Bl.up}}_j \circ T^{-1}_{j, k_{s + 1}})(\xi_{k_{s + 1}}, \eta_{k_{s + 1}})$ for some $j \in \{ k_s + 1, \dots, k_{s + 1} - 1 \}$ preserves the form of $\xi_{k_{s + 1}}$ and $\eta_{k_{s + 1}}$. Abbreviating as in (\ref{Transition}) by $m := -l > 0$, $a := a_j$, $b := b_j$, $(w', z') := (T_{j, k_{s + 1}} \circ T^{\text{Bl.up}}_j \circ T^{-1}_{j, k_{s + 1}})(w_{k_{s + 1}}, z_{k_{s + 1}})$, $(\xi', \eta') := (T_{j, k_{s + 1}} \circ T^{\text{Bl.up}}_j \circ T^{-1}_{j, k_{s + 1}})(\xi_{k_{s + 1}}, \eta_{k_{s + 1}})$, we obtain 

\begin{eqnarray*}
\xi' &=& \xi_i(\xi^q_i\eta^m_i - c)^m \\
&=& (w_i + z_i^{t - s - 1}R) \cdot [(w_i + z_i^{t - s - 1}R)^qz_i^m\underbrace{(1 + z_i^{t - s - 1}S)^m}_{=: 1 + z_i^{t - s - 1}S'(w_i, z_i)} - c ]^m \\
&=& \left( \frac{w'}{(w'^qz'^m - c)^m} + z'^{t - s - 1}\underbrace{(w'^qz'^m - c)^qR\left( \frac{w'}{({w'}^q{z'}^m - c)^m}, z'({w'}^q{z'}^m - c)^q \right)}_{=: \tilde{R}(w', z')} \right) \\
&\cdot& \left[ \left( \frac{w'}{(w'^qz'^m - c)^m} + z'^{t - s - 1}(w'^qz'^m - c)^qR( \cdots ) \right)^q \cdot z'^m(w'^qz'^m - c)^{qm}(1 + z'^{t - s - 1}\underbrace{(w'^qz'^m - c)^{q(t - s - 1)}S'}_{=: S''(w', z')}) - c \right]^m \\
&=& \left( \frac{w'}{(w'^qz'^m - c)^m} + z'^{t - s - 1}\tilde{R} \right) \cdot [ ( w' + z'^{t - s - 1}(w'^qz'^m - c)^{q + m}R )^q \cdot z'^m(1 + z'^{t - s - 1}S'') - c ]^m  
\end{eqnarray*} 

\noindent Using binomial expansion we see that the expression inside the second bracket is of the form $w'^qz'^m - c + z'^{t - s - 1}T$ for some power series $T$. Thus,

\begin{eqnarray*}
\xi' &=& \left( \frac{w'}{(w'^qz'^m - c)^m} + z'^{t - s - 1}\tilde{R} \right) \cdot [ w'^qz'^m  - c + z'^{t - s - 1}T ]^m \\
&=& w' + z'^{t - s - 1}R'
\end{eqnarray*} 

\noindent for some power series $R'$. Repeating this for all $j \in \{ k_s + 1, \dots, k_{s + 1} - 1 \}$, the coordinates $(w_{k_{s + 1}}, z_{k_{s + 1}})$ become $(x_{k_{s + 1}}, y_{k_{s + 1}})$ and the claim follows for $\xi_{k_{s + 1}}$. Similarly we proceed for $\eta_{k_{s + 1}}$. It remains to check that the form of $\xi_{k_{s + 1}}$ and $\eta_{k_{s + 1}}$ does not change after creating a feather on $C_{k_{s + 1}}$, say in $(w_{k_{s + 1}}, z_{k_{s + 1}}) = (c, 0)$. We assume in addition that $s \leq t - 2$ (the case $s = t - 1$ already corresponds to assertion (2) of Claim 1). Since $h(c, 0) = (c, 0)$, we introduce new coordinates via

$$(w', z') = \left( w_{k_{s + 1}}, \frac{z_{k_{s + 1}}}{w_{k_{s + 1}} - c} \right) \quad \text{and} \quad (\xi', \eta') = \left( \xi_{k_{s + 1}}, \frac{\eta_{k_{s + 1}}}{\xi_{k_{s + 1}} - c} \right).$$

\noindent It follows that

$$\xi' = w' + z'^{t - s - 1}(w' - c)^{t - s - 1}R(w', z'(w' - c)) = w' + z'^{t - s - 1}R'(w', z'),$$

\noindent and, using the geometric series,

\begin{eqnarray*}
\eta' &=& \frac{z'(w' - c)(1 + z'^{t - s - 1}(w' - c)^{t - s - 1}S(w', z'(w' - c)))}{w' + z'^{t - s - 1}(w' - c)^{t - s - 1}R(w', z'(w' - c)) - c} \\
&=& \frac{z'(1 + z'^{t - s - 1}(w' - c)^{t - s - 1}S(w', z'(w' - c)))}{1 + z'^{t - s - 1}(w' - c)^{t - s - 2}R(w', z'(w' - c))} \\
&=& z'(1 + z'^{t - s - 1}S'(w', z'))
\end{eqnarray*}

\noindent for some power series $R', S'$. Hence the proof of (1) is complete.
 
Let us show assertion (2). It is shown above that $(\xi_{k_s}, \eta_{k_s})$ in the case $s = t$ has the desired form if no feather is attached to $C_{k_s}$. Let $s \geq t$. By induction hypothesis we have

$$(\xi_{k_s}, \eta_{k_s}) = (x_{k_s} + f_i(a) + y_{k_s}R, y_{k_s}(1 + y_{k_s}S)).$$

\noindent We perform an outer blow-up in some point $(x_{k_s}, y_{k_s}) = (c, 0)$ to create the component $C_{k_{s + 1}}$. Since $h$ maps this point to $(\xi_{k_s}, \eta_{k_s}) = (c + f_i(a), 0)$, we introduce the affine coordinates $(w_{k_{s + 1}}, z_{k_{s + 1}})$ respectively $(\xi_{k_{s + 1}}, \eta_{k_{s + 1}})$ via 

$$(x_{k_s}, y_{k_s}) = (w_{k_{s + 1}}z_{k_{s + 1}} + c, z_{k_{s + 1}}) \quad \text{resp.} \quad (\xi_{k_s}, \eta_{k_s}) = (\xi_{k_{s + 1}}\eta_{k_{s + 1}} + c + f_i(a), \eta_{k_{s + 1}}).$$

\noindent Letting for brevity $(w', z') := (w_{k_{s + 1}}, z_{k_{s + 1}})$, we obtain

\begin{eqnarray*}
(\xi_{k_{s + 1}}, \eta_{k_{s + 1}}) &=& \left( \frac{\xi_{k_s} - c - f_i(a)}{\eta_{k_s}}, \eta_{k_s} \right) \\
&=& \left( \frac{x_{k_s} - c  + y_{k_s}R(x_{k_s}, y_{k_s})}{y_{k_s}(1 + y_{k_s}S(x_{k_s}, y_{k_s}))}, y_{k_s}(1 + y_{k_s}S(x_{k_s}, y_{k_s})) \right) \\
&=& \left( \frac{w'z'  + z'R(w'z' + c, z')}{z'(1 + z'S(w'z' + c, z'))}, z'(1 + z'S(w'z' + c, z')) \right) \\
&=& \left( \frac{w' + R(w'z' + c, z')}{1 + z'S(w'z' + c, z')}, z'(1 + z'S(w'z' + c, z')) \right).
\end{eqnarray*}

\noindent $\eta_{k_{s + 1}}$ already has the desired form. The power series $R(w'z' + c, z') \in \c[[w', z']]$ can be written in the form $R(w'z' + c, z') = b + z'Q(w', z')$ for some $Q \in \c[[w', z']]$. This yields

\begin{eqnarray*}
\xi_{k_{s + 1}} &=& \frac{w' + b + z'Q(w', z')}{1 + z'S(w'z' + c, z')} \\
&=& (w' + b + z'Q(w', z'))(1 + z'S'(w'z' + c, z')) \\
&=& w' + b + z' \cdot [w'S'(w'z' + c, z') + bS'(w'z' + c, z') + Q(w', z') + z'Q(w', z')S'(w'z' + c, z')]
\end{eqnarray*}
 
\noindent where $(1 + z'S(w'z' + c, z'))^{-1} = 1 + z'S'(w'z' + c, z')$. Hence, $\xi_{k_{s + 1}}$ is again of the form $w' + a' + z'R'(w', z')$ with some constant $a'$, depending on $a$. Now we perform inner blow-ups to create the components $C_{k_s + 1}, \dots, C_{k_{s + 1} - 1}$ and then further blow-ups to create the feathers. Again, we have to show that replacing $(w_{k_{s + 1}}, z_{k_{s + 1}})$ by $(T_{j, k_{s + 1}} \circ T^{\text{Bl.up}}_j \circ T^{-1}_{j, k_{s + 1}})(w_{k_{s + 1}}, z_{k_{s + 1}})$ for some $j \in \{ k_s + 1, \dots, k_{s + 1} - 1 \}$ (and similarly for $(\xi_{k_{s + 1}}, \eta_{k_{s + 1}})$) preserves the form of $(\xi_{k_{s + 1}}, \eta_{k_{s + 1}})$ asserted in (2). However, this is a similar computation as in the corresponding part of the proof of assertion (1). Further, it is straightforward to check that the form of $(\xi_{k_s}, \eta_{k_s})$ remains the same if we create additional feathers on $C_{k_s}$ by blow-ups of points on $C_{k_s} \backslash C_{k_s - 1}$. Hence (2) follows.

Finally, it remains to show (3) in the case that $s \geq t$. Again, using assertion (2) it is easy to check that formula (\ref{ShiftInner}) holds as well for $s \geq t$, and even, if we create further feathers on $C_j$ for some $j \in \{ k_s + 1, \dots, i - 1 \}$. Hence we obtain assertion (3) for every $s = 2, \dots, r - 1$. This completes the proof of Claim 1.

Now we show that any elementary shift $h(x_0, y_0) = (x_0 + ay_0^{t - 2}, y_0)$ maps the intersection point $F \cap F^\vee$ on the intersection point $h(F) \cap h(F)^\vee = G \cap G^\vee$, if $F$ is attached to $C_i$, where $C_i$ is inner and not exceptional. We fix an arbitrary $i \in \{ k_s + 1, \dots, k_{s + 1} - 1 \}$ and assume that all feathers attached to any $C_\mu$ with $\mu < i$ are already created. Therefore we use the coordinates $(w_i, z_i)$ on the corresponding intermediate surface. To create the feather $F$ on the component $C_i$ we blow up in a point $(w_i, z_i) = (d, 0) \in C_i \backslash (C_{i - 1} \cup C_{i + 1})$. By (3), the map $h$ fixes the point $(d, 0)$. After the blow-up in $(d, 0)$ we introduce the affine coordinates $(u, v)$ and $(\xi, \eta)$, respectively, via $(w_i, z_i) = (uv + d, v)$ and $(\xi_i, \eta_i) = (\xi\eta + d, \eta)$, respectively. In these coordinates we have $F \backslash D = \{ v = 0 \}$ and $h(F) \backslash h.D = \{ \eta = 0 \}$. A short direct computation shows that the lift of $h$ has the form

\begin{eqnarray*}
(\xi, \eta) &=& (u + v^{b_i - 1}A(u, v), v \cdot (1 + v^{b_i}B(u, v)))
\end{eqnarray*}
 
\noindent for certain power series $A$ and $B$ (and furthermore, $A(u, 0) = d^{a_i + 1} \cdot R(d, 0)$). Thus the map $h \vert_{F \backslash D}: F \backslash D \to h(F) \backslash h.D$ is given by

\begin{equation}\label{AutOnFeather}
h(u, 0) = \begin{cases} (u + A(u, 0), 0) &, b_i = 1 \\ \left(u, 0 \right) &, b_i \geq 2.\end{cases}
\end{equation}

\noindent Hence, if $F$ is not attached to an exceptional component, the map $h$ induces the identity on $F$ (recall that $b_i = 1$ if and only if $C_i$ is exceptional). In particular, $h(0, 0) = (0, 0)$ holds for every elementary shift $h$. But since $(u, v) = (0, 0)$ equals $F \cap F^\vee$ on $X$ and $(\xi, \eta) = (0, 0)$ equals $h(F) \cap h(F^\vee)$ on $(X', D')$, the claim follows immediately.\\
\end{bew}

Another technical subtlety (which shows that feathers attached either to exceptional components or to outer components do not provide points in the complement of $O$), which plays an important role in the proof of the main theorem (see Theorem \ref{MainTheorem}) is stated in the following lemma:

\begin{lemma}\label{TranslationOnFeather} Let $(X, D)$ be a $(-1)$-completion of a smooth Gizatullin surface $V$, $F$ a feather of $D_{\text{ext}}$ which is attached to $C_i$ and let $h$ be the elementary shift $h_{a, t + 1}(x_0, y_0) = (x_0 + ay_0^{t - 1}, y_0)$. If $i = k_t$ for some $t$ or if $i \in \mathfrak{E}_D$ with $k_t < i < k_{t + 1}$, then the lift of $h$,

$$h: F \backslash D \to h(F) \backslash h.D$$

\noindent has the form $h(u, 0) = (u + f(a), 0)$ (in the coordinates $(u, v)$ and $(\xi, \eta)$ used in (\ref{AutOnFeather})), where $f(a)$ is non-zero constant for general $a$.\\
\end{lemma}

\begin{bew} \textit{Case 1: $C_i$ is an outer component, $i = k_t$ for some $t = 2, \dots, r$.} We consider the lift of $h$ onto an intermediate surface $X_M$, $M \leq N$, where $X_M$ contains only feathers with mother component $C_j$, $j \leq k_t - 1$. Then by Claim 1 in the proof of Lemma \ref{IntersectionPointOfFeathersByShifts}, $h$ has the form

$$(\xi_{k_t}, \eta_{k_t}) = h(w_{k_t}, z_{k_t}) = (w_{k_t} + z_{k_t}R(w_{k_t}, z_{k_t}), z_{k_t}(1 + z_{k_t}S(w_{k_t}, z_{k_t}))).$$

\noindent The same computation which leads to (\ref{AutOnFeather}) shows, that the claim is equivalent to the following property: $R(u, 0)$ is a non-trivial polynomial in $a$ for every fixed $u \in \c$.

First, this is obviously true if $X = Q$ (and $t = 3$). Now, considering the inductive contruction of our local coordinate charts  and the computations made above it is easy to see that this property is preserved after every step. Hence the claim follows.\\

\noindent \textit{Case 2: $C_i$ is an inner exceptional component.} The claim follows by a similar argument as in case 1.
\end{bew}

For our main result we will have to assume that $V$ admits a $(-1)$-completion $(X, D)$ such that its reversion $(X^\vee, D^\vee)$ is also a $(-1)$-completion (be careful, a reversion of a $(-1)$-completion does not necessarily result in a $(-1)$-completion!). Fortunately, we get the second condition for free, assuming only that $V$ admits a $(-1)$-completion.

\begin{lemma} (cf. also \cite{FKZ4}, Prop. 4.3.8 for a special case). Let 

$$X_n = X(M_2, \dots, M_{k_3 - 1}, c_{k_3}, M_{k_3}, \dots, c_n, M_n)$$

\noindent be a $(-1)$-presentation. Then there exists a finite sequence of elementary shifts which transforms $X_n$ into a new $(-1)$-presentation 

$$\tilde{X}_n = X(\tilde{M}_2, \dots, \tilde{M}_{k_3 - 1}, \tilde{c}_{k_3}, \tilde{M}_{k_3}, \dots, \tilde{c}_n, \tilde{M}_n),$$

\noindent such that its reversion is again of $(-1)$-type.
\end{lemma}

\begin{bew} In essence, we follow the proof of Prop. 4.3.8 in \cite{FKZ4}, which asserts the same for the special case where all components of $D$ are outer ones.\\

\noindent \textbf{Claim 1:} Let $X_n = X(M_2, \dots, M_{k_3 - 1}, c_{k_3}, M_{k_3}, \dots, c_n, M_n)$ be a $(-1)$-presentation and $h_{a, t}(x_0, y_0) = (x_0 + ay_0^{t - 2}, y_0)$. Then $h.X_n$ has the presentation

$$h.X_n = X(M_2, \dots, M_{k_t - 1}, c_{k_t}, a + M_{k_t}, M'_{k_t + 1} \dots, M'_{k_{t + 1} - 1}, a + c_{k_{t + 1}}, M'_{k_{t + 1} + 1}, \dots, c'_n, M'_n)$$

\noindent for some $c'_{k_i}$ and $M'_j$ for $i > t + 1$ and $j > k_t$.\\

\noindent \textit{Proof of Claim 1:} For $t = 2$ the claim is obviously true. Hence we assume in the sequel that $t \geq 3$. Since $X_n$ is a $(-1)$-presentation, it can be obtained by first creating the sub-zigzag $C_0 \cup C_1 \cup C_2 \cup C_{k_3} \cup \cdots \cup C_{k_{r - 1}} \cup C_n$ (which consists only of the outer components), then by blowing up all $M_{k_i} \subseteq C_{k_i} \backslash C_{k_{i - 1}}$, $2 \leq i \leq r$, and finally by creating all inner components and the corresponding feathers. We consider again our coordinates $(x_0, y_0)$ on $\a^2 = Q \backslash (C_0 \cup C_1)$. After a suitable translation we may suppose that $c_{k_3} = (0, 0)$. The blow-up with center $c_{k_3}$ can be written in coordinates as

$$(x', y') = (x_0/y_0, y_0).$$

\noindent In these coordinates the exceptional curve $C_{k_3}$ is given by $y' = 0$ and the proper transform of $C_2$ by $x' = \infty$. Hence, $h_{a, t}$ lifts as

\begin{equation}\label{LiftOfShift}
h_{a, t}: (x', y') \mapsto (x' + a{y'}^{t - 3}, y').
\end{equation}

\noindent In particular, $h_{a, 3}$ yields the identity on $C_{k_3} \backslash C_2 \cong \a^1$, and consequently, on all $C_3, \dots, C_{k_3 - 1}$. Otherwise, if $t = 3$, it yields the translation by $a$ on $C_{k_3} \backslash C_2 \cong \a^1$, and the identity on all $C_3, \dots, C_{k_3 - 1}$. Formula (\ref{LiftOfShift}) remains the same after replacing the coordinates $(x', y')$ by the new ones  $(x' - c, y')$, where $c_{k_4} = (c, 0)$. Thus we may assume that $c_{k_4} = (0, 0)$ in the coordinate system $(x', y')$. Now the claim follows by induction.\\  

\noindent \textbf{Claim 2:} Let 

$$X_n = X(M_2, \dots, M_{k_t - 1}, c_{k_t}, M_{k_t}, \dots, c_n, M_n)$$

\noindent and 

$$X'_n = X(M_2, \dots, M_{k_t - 1}, c_{k_t}, M'_{k_t}, \dots, c'_n, M'_n)$$

\noindent be presentations with reversed presentation

$$X^\vee_n = X(M_n, \dots, M_{k_t + 1}, c^\vee_{k^\vee_t}, M_{k_t}, \dots, c^\vee_n, M_2)$$

\noindent and 

$${X'}^\vee_n = X(M'_n, \dots, M'_{k_t + 1}, {c'}^\vee_{k^\vee_t}, M'_{k_t}, M_{k_t - 1} \dots, {c'}^\vee_n, M_2).$$

\noindent Then we have

$$c^\vee_{k^\vee_i} = {c'}^\vee_{k^\vee_i} \quad \forall \ i = t + 1, \dots, r - 1.$$

\noindent \textit{Proof of Claim 2:} Starting with the completion $(X_n, D_n)$ of $V$, we consider the correspondence fibration $\psi: (W, E) \to \p^1$ for the pair $(C_t, C^\vee_{t^\vee})$. Similarly, we let $\psi': (W', E') \to \p^1$ be the correspondence fibration associated to $(X'_n, D'_n)$ for the pair $(C'_t, {C'}^\vee_{t^\vee})$. To obtain the part ${D^\vee}^{\geq t^\vee} = C^\vee_n \cup \cdots \cup C^\vee_{t^\vee}$ of the reversed zigzag $D^\vee$, only inner elementary transformations with centers at the components $C_0 = C'_0$, ..., $C_{k_t - 1} = C'_{k_t - 1}$ are required. Hence, it follows that $C^\vee_{i^\vee} = {C'}^\vee_{i^\vee}$ for all $i^\vee \geq k^\vee_t + 1$. In particular, Claim 2 holds.\\

Now, the assertion of the lemma can be shown as follows. Let 

$$X^\vee_n = X(M_n, \dots, M_{k^\vee_{r - 1} + 1}, c^\vee_{k^\vee_{r - 1}}, M_{k^\vee_{r - 1}}, \dots, c^\vee_n, M_2)$$

\noindent be the reversion of $X_n$. With a suitable coordinate on $\a^1 \cong C_{k_t} \backslash C_{k_t - 1}$, the elementary shift $h_{a, t}$ transforms $X_n$ into a $(-1)$ presentation

$$X'_n = X(M_2, \dots, M_{k_t - 1}, c_{k_t}, a + M_{k_t}, M'_{k_t + 1}, \dots, M'_{k_{t + 1} - 1}, a + c_{k_{t + 1}}, M'_{k_{t + 1}}, \dots, c'_n, M'_n)$$

\noindent with reversion 

$${X'}^\vee_n = X(M'_n, \dots, M'_{k_t - 1}, {c'}^\vee_{k^\vee_t}, a + M_{k_t}, M_{k_t - 1}, \dots, M_{k_{t - 1} + 1}, c^\vee_{k^\vee_{t - 1}}, M_{k_{t - 1}}, \dots, c^\vee_n, M_n),$$

\noindent see Claim 2. Choosing $a$ general we may suppose that

\begin{equation}\label{NewPresentation}
c^\vee_{k^\vee_{t - 1}} \not\in \tilde{M}_t := a + M_t.
\end{equation}

\noindent Applying successively shifts $h_{a_t, t}$, $t = 3, \dots, r$, with general $a_t \in \c$, the resulting surface $\tilde{X}_n$ satisfies (\ref{NewPresentation}) for all $t = 3, \dots, r$, hence it is of $(-1)$-type.
\end{bew}

\begin{rem} As mentioned above, reversing a $(-1)$-completion $(X, D)$ does not result in a $(-1)$-completion in general. Thus, studying the automorphisms of smooth Gizatullin surfaces which admit a $(-1)$-completion forces us to consider also lifts of elementary shifts on completions of $V$ which are not $(-1)$-completions. However, a similar computation as in the proof of Lemma \ref{IntersectionPointOfFeathersByShifts} yields that $\xi_i$ and $\eta_i$ are of the same form if $k_s < i < k_{s + 1}$, disregarding the fact that some of the feathers may be \emph{not} $(-1)$-feathers. This is all we need to prove our main result below.
\end{rem}

\bigskip
In the following we fix some notations. Let $W_m := \{ z \in \c^* \mid z^m = 1 \} \cong \z_m$ be the set (group) of $m$-th roots of unity. Every finite non-empty subset $A \subseteq \c^*$ can be written in the form $A = \bigcup_{i = 1}^s c_iW_{m_i}$, $c_i \in \c^*$, such that $s$ is minimal. We let

$$G(A) := \{ \alpha \in \c^* \mid \alpha \cdot A = A \} \cong \z_{\gcd(m_1, \dots, m_s)} \quad \text{and} \quad m(A) := \frac{m_1 + \cdots + m_s}{\gcd(m_1, \dots, m_s)}.$$

The integer $m(A)$ is precisely the number of $G(A)$-orbits under the $G(A)$-action on $A$. In addition, if $A = \emptyset$, we let

$$G(\emptyset) := \c^* \quad \text{and} \quad m(\emptyset) := 0.$$

Using Lemma \ref{IntersectionPointOfFeathersByShifts} we can exhibit further examples of smooth quasi-homogeneous surfaces which are not homogeneous. The following theorem is the main result of this article:

\begin{thm}\label{MainTheorem} Let $V$ be a smooth Gizatullin surface which admits a $(-1)$-completion and let $(X, D)$ a standard completion of $V$ of $(-1)$-type. Let $A_i = \{ P_{i, 1}, \dots, P_{i, r_i} \} \subseteq C_i \backslash (C_{i - 1} \cup C_{i + 1}) \cong \c^*$, $i \not\in \{ k_2, k_3, \dots, k_r \}$, be the base point set of the feathers $F_{i, j}$. Furthermore, for every $i$ with $k_s + 1 \leq i \leq k_{s + 1} - 1$ and $i \not\in \mathfrak{E}_D \cup \mathfrak{E}^\vee_{D^\vee}$ we let $B_{i, 1}, \dots, B_{i, m(A_i)}$ be the orbits of the $G(A_i)$-action on $A_i$. If the configuration invariant $Q(X, D) = (Q_2, \dots, Q_n)$ of $V$ is not symmetric (\ie $Q_i \neq Q_{i^\vee}$ for some $i$), then we let 

\begin{equation}\label{InvariantSubsets1}
O_{i, j} := \bigcup_{1 \leq l \leq r_i; P_{i, l} \in B_{i, j}} (F_{i, l} \cap F^\vee_{i, l}) \subseteq V, \quad 1 \leq j \leq m(A_i).
\end{equation}

\noindent Otherwise, we let for $i \leq \lfloor \frac{n}{2} \rfloor + 1$

\begin{equation}\label{InvariantSubsets2}
O_{i, j} := \bigcup_{1 \leq l \leq r_i; P_{i, l} \in B_{i, j}} (F_{i, l} \cap F^\vee_{i, l}) \cup \bigcup_{1 \leq l \leq r_i; P_{i^\vee, l} \in B_{i^\vee, j}} (F_{i^\vee, l} \cap F^\vee_{i^\vee, l}) \subseteq V, \quad 1 \leq j \leq m(A_i),
\end{equation}

\noindent where we identify $B_{i, j}$ with $B_{i^\vee, j}$ after a suitable numbering of the orbits.\footnote[2]{Since $Q_i = Q_{i^\vee}$, we have $A_{i^\vee} = \alpha_i \cdot A_i$ for some $\alpha_i \in \c^*$. Thus we have $B_{i^\vee, j} = \alpha_i \cdot B_{i, j}$ for a suitable numbering of the orbits. Note, that this correspondence is well-defined, since two such $\alpha_i$ differ by an element of $G(A_i) = G(A_{i^\vee})$, which leave the $B_{i^\vee, j}$ invariant.} Moreover, in both cases, we let

$$O_0 := V \backslash \left( \bigcup_{i, j} O_{i, j} \right).$$

\noindent Then the following hold:
\begin{itemize}
\item[(1)] The subsets $O_0$ and $O_{i, j}$ are invariant under the action of $\Aut(V)$. Moreover, $O_0$ contains the big orbit $O$.
\item[(2)] Let $F$ be a feather of $D_{\text{ext}}$ which is attached to $C_i$, such that either $C_i$ is an outer component or $i \in \mathfrak{E}_D \cup \mathfrak{E}^\vee_{D^\vee}$ holds. Then $F \backslash D$ is contained in $O$. 
\item[(3)] Assume, that $r_i > 0$ holds for a unique $i \in \{ 2, \dots, n \}$. Then the subsets $O_0$ and $O_{i, j}$ form the orbit decomposition of the natural action of $\Aut(V)$ on $V$. In particular, $O = O_0$ holds.
\end{itemize}
\end{thm}

\begin{bew} Let $\varphi \in \Aut(V)$. We choose an arbitrary standard completion $(X, D)$ of $V$ and extend $\varphi$ to a birational map

$$\bar{\varphi}: (X, D) \dasharrow (X, D).$$

\noindent We decompose $\bar{\varphi}$ into a sequence $\bar{\varphi} = \varphi_m \circ \cdots \circ \varphi_1$ according to Prop. \ref{FactorizationOfBirationalMapsStandard}. In each of these completions $(X_k, D_k)$, $0 \leq k \leq m$, we denote by $O^{(k)}_{i, j}$ the corresponding subsets defined in (\ref{InvariantSubsets1}) or (\ref{InvariantSubsets2}) (and similarly for the boundary components, the base point sets of the feathers etc.). Now there are three possible cases:\\
\noindent \textbf{Case 1: $\bf\varphi_k$ is a reversion.} Let $\varphi_k: (X_{k - 1}, D_{k - 1}) \dasharrow (X_k, D_k)$ be a reversion. After a suitable numbering of the feathers in $(X_k, D_k)$, the proper transforms of $F^{(k - 1)}_{i, j}$ and $\left( {F_{i, j}^{(k - 1)}} \right)^\vee$, respectively, under $\varphi_k$ are $\left( {F_{i^\vee, j}^{(k)}} \right)^\vee$ and $F^{(k)}_{i^\vee, j}$, respectively (take a minimal resolution $(X_{k - 1}, D_{k - 1}) \leftarrow (Z, B) \to (X_k, D_k)$ of $\varphi_k$). Hence, $\varphi_k$ maps (after a suitable numbering) $F_{i, j}^{(k - 1)} \cap \left( {F_{i, j}^{(k - 1)}} \right)^\vee$ onto $F_{i, j}^{(k)} \cap \left( {F_{i, j}^{(k)}} \right)^\vee$, which implies that after a suitable numbering of the $G(A_i^{(k)})$-orbits we have $\varphi_k(O^{(k - 1)}_{i, j}) = O^{(k)}_{i^\vee, j}$.\\
\noindent \textbf{Case 2: $\bf\varphi_k$ is an elementary shift.} By Lemma \ref{IntersectionPointOfFeathersByShifts}, the map $\varphi_k$ maps $O^{(k - 1)}_{i, j}$ onto $O^{(k)}_{i, j}$.\\
\noindent \textbf{Case 3: $\bf\varphi_k$ is a torus element.} First, note that fibered maps do not change the configuration invariant $Q(V, D)$. In particular, $\varphi_k$ preserves $Q(V, D)$.

Let $\varphi_k(x_0, y_0) = (ax_0, by_0)$ for some $a, b \in \c^*$. Introducing the same coordinates $(x_i, y_i)$, $(w_i, z_i)$ and their images $(\xi_i, \eta_i)$ and observing the behaviour of $\varphi_k$ in these coordinates it is easy to check that $\varphi_k$ induces on \emph{every} boundary component a multiplication $u \mapsto \lambda u$ for some constant $\lambda$ (depending on $a, b$ and the component). We fix some $i \in \{ k_s + 1, \dots, k_{s + 1} - 1 \}$. The map $\varphi_k: C^{(k - 1)}_i \backslash (C^{(k - 1)}_{i - 1} \cup C^{(k - 1)}_{i + 1}) \cong \c^* \to C^{(k)}_i \backslash (C^{(k)}_{i - 1} \cup C^{(k)}_{i + 1}) \cong \c^*$ has the form $w^{(k - 1)}_i \mapsto \lambda w^{(k)}_i$. Since $\varphi_k$ maps the base point sets $A^{(k - 1)}_i$ onto $A^{(k)}_i$ and since the configuration invariants $Q(X_{k - 1}, D_{k - 1})$ and $Q(X_k, D_k)$ are equal, two such $\lambda$ differ only by an element of $G(A^{(k - 1)}_i)$ (which is by the Matching Principle either equal to $G(A_i)$ or to $G(A_{i^\vee})$). Therefore, identifying $A^{(k - 1)}_i$ with $A^{(k)}_i$, the base points of the feathers $F^{(k - 1)}_{i, j}$ and $\varphi_k(F^{(k - 1)}_{i, j})$ are in the same $G(A^{(k - 1)}_i)$-orbit. Now, lifting $\varphi_k$ to the affine part of the feathers, $\varphi_k\vert_{F \backslash D}: F \backslash D \to \varphi_k(F) \backslash \varphi_k.D$ has the form $(u, 0) \mapsto (\lambda' u, 0)$ for some $\lambda' \in \c^*$, and hence $\varphi_k(0, 0) = (0, 0)$. In other words, 

$$\varphi_k(F^{(k - 1)}_{i, j} \cap (F^{(k - 1)}_{i, j})^\vee) = \varphi_k(F^{(k - 1)}_{i, j}) \cap (\varphi_k(F^{(k - 1)}_{i, j}))^\vee.$$

\noindent It follows that $\varphi_k$ maps $O^{(k - 1)}_{i, j}$ onto $O^{(k)}_{i, \sigma(j)}$. As in case 1 we may assume that $\varphi_k(O^{(k - 1)}_{i, j}) = O^{(k)}_{i, j}$.

We consider now the configuration invariant $Q(X, D) = (Q_2, \dots, Q_n)$ of $V$. Here we have to distinguish two cases, namely either $Q(X, D)$ is symmetric (\ie $Q_i = Q_{i^\vee}$ for all $i$) or not. Assume first, that $Q(X, D)$ is not symmetric. Then the number of reversions among $\varphi_1, \dots, \varphi_m$ is even and the above cases imply that $\bar{\varphi} = \varphi_m \circ \cdots \circ \varphi_1$ maps $O_{i, j} = O^{(0)}_{i, j}$ onto $O^{(k)}_{i, j} = O_{i, j}$. Hence, the subsets $O_{i, j}$ are invariant under $\bar{\varphi}\vert_V = \varphi$. If $Q(X, D)$ is symmetric, then the number of reversions occuring in the decomposition of $\bar{\varphi}$ may also be odd. However, by construction of the $O_{i, j}$ in the symmetric case we have $\varphi_k(O_{i, j}^{(k - 1)}) = O_{i, j}^{(k)}$ (again, after a suitable numbering of the $G(A_i^{(k)})$-orbits). From here we conclude as above. Hence, assertion (1) follows.

To show (2), we consider a feather $F$ which is either attached to an outer component $C_{k_{s - 1}}$ or to an inner exceptional component $C_i$ with $k_{s - 1} < i < k_s$, and a point $p \in F \backslash D \subseteq V$. By Lemma \ref{TranslationOnFeather}, the elementary shift $h(x_0, y_0) = h_{a, s}(x_0, y_0) = (x_0 + ay^{s - 2}_0, y_0)$, $a \neq 0$, induces a non-trivial translation on $F \backslash D$ (in 
the coordinates $(u, v)$ and $(\xi, \eta)$ used in (\ref{AutOnFeather}), see Lemma \ref{TranslationOnFeather}). In other words, $h$ induces for general $a \in \c$ a birational map

$$h: (X, D) \dasharrow (X', D') = h.(X, D),$$

\noindent which induces an isomorphism $V = X \backslash D \stackrel{\sim}{\to} V' = X' \backslash D'$ by restriction, such that $h(p)$ is contained in $h(F)$, but not in any matching feather. By Prop. \ref{AutFiniteComplement}, $h(p)$ is contained in the big orbit of $V'$ and therefore, the same holds for $p$. Hence (2) follows.

Finally, assume that $r_i > 0$ holds for a unique $i \in \{ 2, \dots, n \}$. If $n \geq 4$ , then $C_3, \dots, C_{n - 1}$ are necessarily $*$-components and all feathers are attached to some $C_i$ with $3 \leq i \leq n - 1$. In this case assertion (4) is precisely the statement of \cite{Ko}, Theorem 3.11 (4). Otherwise, if $n \leq 3$, then we necessarily have either $n = 1$ (which gives $V = \a^2$) or $n = 2$ and all feathers are attached to $C_2$. In the last case we obtain the Danielewski surfaces $V = \{ xy - P(z) = 0 \} \subseteq \a^3$, where $P(z)$ is a polynomial with pairwise distinct roots (see \cite{BD1}, Thm. 5.4.5). However, in both cases the action of $\Aut(V)$ is transitive, \ie $O = V = O_0$.
\end{bew}

Now the following corollary is obvious:

\begin{cor}\label{CorClassification} Let $V$ be a smooth Gizatullin surface that admits a $(-1)$-completion. Consider a presentation $X_n = X(M_2, M_3, \dots, c_{k_3}, M_{k_3}, \dots, c_n, M_n)$ of $V$. Then $V$ is not homogeneous if there exists an integer $i \in \{2, \dots, n\}$ with $i \not\in \mathfrak{E}_D \cup \mathfrak{E}^\vee_{D^\vee} \cup \{ k_2, k_3, \dots, k_r \}$ such that there is at least one feather attached to $C_i$.
\end{cor}

\begin{rem} It is much more involved to give a necessary condition for $V$ to be not homogeneous. The reason is as follows. By Prop. \ref{AutFiniteComplement}, it is still possible that certain points $F_{i, \rho} \cap F^\vee_{j, \sigma}$ with $i \neq j$ are contained in $V \backslash O$. Indeed, if this intersection is not empty and if $C_i, C_j, C^\vee_{i^\vee}$ and $C^\vee_{j^\vee}$ are neither outer components, nor exceptional ones, then the proof of Lemma \ref{IntersectionPointOfFeathersByShifts} shows that $F_{i, \rho} \cap F^\vee_{j, \sigma}$ is indeed contained in the complement of $O$ (since all elementary shifts induce the identity on the affine part of the $F_{i, \rho}$ and $F^\vee_{j, \sigma}$). However, it is still not clear whether these feathers do intersect or not. At least, the algorithm given in section \ref{PresentationGizatullinSurface} allows us to compute explicit equations for all feathers and their duals, hence it is possible to decide this for any concrete presentation of $V$. 
\end{rem}

\section{Automorphism group of Gizatullin surfaces with a Rigid Extended Divisor}\label{SectionHugeAut}

This section is devoted to the study of the structure of automorphism groups of smooth Gizatullin surfaces. We exhibit among non-homogeneous smooth Gizatullin surfaces examples of surfaces with a huge group of automorphisms. Hugeness means that the automorphism group has the following properties: (1) the subgroup $\Aut(V)_{\text{alg}}$, generated by all algebraic subgroups, is not generated by a countable subset of algebraic subgroups and (2) the quotient $\Aut(V)/\Aut(V)_{\text{alg}}$ contains a free group over an uncountable set of generators (see also \cite{BD2}, 4.2).

To exhibit these desired properties for certain Gizatullin surfaces $V$, it is more convenient to consider $1$-standard completions of $V$, which arise by a blow-up of the Hirzebruch surface $\f_1$. In the following we introduce some notations concerning $1$-standard pairs. According to \cite{BD1}, Lemma 1.0.7, any smooth $1$-standard pair $(X, D)$ may be obtained by some blow-ups of points on a fiber of $\f_1$. Let us fix an embedding of $\f_1$ into $\p^2 \times \p^1$ via

$$\f_1 = \{ ((x: y: z), (s: t)) \in \p^2 \times \p^1 \mid yt - zs = 0 \} = \text{Bl}_{(1: 0: 0)}(\p^2).$$

\noindent We denote by $\tau : \f_1 \to \p^2$ the projection on $\p^2$, by $C_0$ and $C_2$ the lines $\{ z = 0 \}$ and $\{ y = 0 \}$ as well as their proper transforms on $\f_1$ and by $C_1$ the exceptional curve $\tau^{-1}(1: 0: 0) = \{(1: 0: 0)\} \times \p^1$. Moreover, we have isomorphisms

$$\a^2 \stackrel{\cong}{\to} \f_1 \backslash (C_0 \cup C_1),\quad (u_0, v_0) \mapsto ((u_0: v_0: 1), (v_0: 1))$$

\noindent as well as 

$$\a^2 \stackrel{\cong}{\to} \p^2 \backslash C_0,\quad (u_0, v_0) \mapsto (u_0: v_0: 1).$$

\noindent In these coordinates we have $C_2 = \overline{\{ v_0 = 0 \}}$. In the following we denote these coordinates on $\f_1 \backslash (C_0 \cup C_1)$ by $(u_0, v_0)$. The $\p^1$-fibration $\rho: \f_1 \to \p^1$, induced by the linear pencil $|C_0|$ is simply given by $(u_0, v_0) \mapsto v_0$ on the affine part $\a^2 = \f_1 \backslash (C_0 \cup C_1)$.\\

We denote by $\Aff$ the group of automorphisms of $\a^2$, which extend to automorphisms of $\p^2$ and by $\Jon$ the group of triangular (de Jonqui\`{e}res) automorphisms, \ie automorphisms which preserve the fibration given by $|C_0|$. In other words, we have

\begin{eqnarray*}
\Aff &=& \{ (u_0, v_0) \mapsto (a_{11}u_0 + a_{12}v_0 + b_1, a_{21}u_0 + a_{22}v_0 + b_2) \mid a_{11}a_{22} - a_{12}a_{21} \neq 0 \} \\
\Jon &=& \{ (u_0, v_0) \mapsto (au_0 + P(v_0), bv_0 + c) \mid a, b \in \c^*, c \in \c, P(v_0) \in \c[v_0] \}.
\end{eqnarray*}

\noindent Sometimes we need the action of the $2$-torus 

$$\t := \{ (u_0, v_0) \mapsto (au_0, bv_0) \mid a, b \in \c^* \} \subseteq \Jon \cap \Aff.$$

\bigskip
Moreover, if we consider a reversion $(X, D) \dasharrow (X', D')$ of $1$-standard pairs centred in a point $p \in C_0 \backslash C_1$, we associate to this point its image $(\lambda: 1: 0)$ in $\p^2$ via the map $\tau \circ \eta: X \to \p^2$. Using this identification we simply write $\lambda \in C_0 \backslash C_1$.

\bigskip
We recall the following usefull lemma:

\begin{lemma}\label{JonquieresAutomorphisms2} (\cite{BD1}, Lemma 5.2.1) For $i = 1, 2$, let $(X_i, D_i, \overline{\pi}_i)$ be a $1$-standard pair with a minimal resolution of singularities $\mu_i: (Y_i, D_i, \overline{\pi}_i \circ \mu_i) \to (X_i, D_i, \overline{\pi}_i)$ and let $\eta_i: Y_i \to \f_1$ be the (unique) birational morphism. Then, the following statements are equivalent:
\begin{itemize}
\item[(a)] The $\a^1$-fibered surfaces $(X_1 \backslash D_1, \pi_1)$ and $(X_2 \backslash D_2, \pi_2)$ (respectively the pairs $(X_1, D_1, \overline{\pi}_1)$ and $(X_2, D_2, \overline{\pi}_2)$) are isomorphic.
\item[(b)] There exists an element of $\Jon$ (respectively of $\Jon \cap \Aff$) which sends the points blown-up by $\eta_1$ onto those blown-up by $\eta_2$ and sends the curves contracted by $\mu_1$ onto those contracted by $\mu_2$.
\end{itemize}
\end{lemma}

Following \cite{BD1}, we introduce for an $\a^1$-fibered surface $V$ a (not necessarily finite) graph $\F_V$, which reflects the structure of the automorphism group of $V$, as follows.

\begin{df} (\cite{BD1}, Def. 4.0.5) To every normal quasi-projective surface $V$ we associate the oriented graph $\F_V$ as follows:
\begin{itemize}
\item[(1)] A vertex of $\F_V$ is an equivalence class of a $1$-standard pair $(X, D)$, such that $X \backslash D \cong V$, where two $1$-standard pairs $(X_1, D_1, \bar{\pi}_1)$ and $(X_2, D_2, \bar{\pi}_2)$ define the same vertex if and only if $(X_1 \backslash D_1, \pi_1) \cong (X_2 \backslash D_2, \pi_2)$.
\item[(2)] An arrow of $\F_V$ is an equivalence class of reversions. If $\varphi: (X, D) \to (X', D')$ is a reversion, then the class $[\varphi]$ of $\varphi$ is an arrow starting from $[(X, D)]$ and ending at $[(X', D')]$. Two reversions $\varphi_1: (X_1, D_1) \dasharrow (X'_1, D'_1)$ and $\varphi_2: (X_2, D_2) \dasharrow (X'_2, D'_2)$ define the same arrow if and only if there exist isomorphisms $\theta: (X_1, D_1) \to (X_2, D_2)$ and $\theta': (X'_1, D'_1) \to (X'_2, D'_2)$, such that $\varphi_2 \circ \theta = \theta' \circ \varphi_1$.
\end{itemize}
\end{df}

\begin{rem} It follows from the definition that for a given $1$-standard pair $(X, D)$, two reversions $\varphi: (X, D) \dasharrow (X_1, D_1)$ and $\varphi: (X, D) \dasharrow (X_2, D_2)$ with centers in $p_1$ and $p_2$, respectively, define the same arrow if and only if there exists an automorphism $\psi \in \Aut(X, D)$ such that $\psi(p_1) = p_2$.
\end{rem}

The structure of the graph $\F_V$ allows us to decide, whether the automorphism group $\Aut(V)$ of $V$ is generated by automorphisms of $\a^1$-fibrations. Here we say that $\varphi \in \Aut(V)$ is an \textit{automorphism of $\a^1$-fibrations} if there exists an $\a^1$-fibration $\pi: V \to \a^1$ and an automorphism $\psi \in \Aut(\a^1)$, such that $\pi \circ \varphi = \psi \circ \pi$. More generally, we say that two $\a^1$-fibered surfaces $(V, \pi)$ and $(V', \pi')$ are isomorphic if there exist isomorphisms $\varphi: V \to V'$ and $\psi: \a^1 \to \a^1$ such that $\pi' \circ \varphi  = \psi \circ \pi$ and write $(V, \pi) \cong (V', \pi')$.

\begin{prop}\label{SurfaceTree}(\cite{BD1}, Prop. 4.0.7) Let $V$ be a normal quasi-projective surface with a non-empty graph $\F_V$. Then the following holds:
\begin{itemize}
\item[(1)] The graph $\F_V$ is connected.
\item[(2)] There is a natural bijection between the set of vertices of $\F_V$ and the isomorphism classes of $\a^1$-fibrations on $V$. 
\item[(3)] Let $(X, D)$ be a $1$-standard pair with $X \backslash D \cong V$ and let $D$ contain at least one curve with self-intersection $\leq -3$. Then there is a natural exact sequence

$$1 \to H \to \Aut(V) \to \Pi_1(\F_V) \to 1,$$

\noindent where $H$ is the (normal) subgroup of $\Aut(V)$ generated by all automorphisms of $\a^1$-fibrations and $\Pi_1(\F_V)$ is the fundamental group of the graph $\F_V$. In particular, the graph $\F_V$ is a tree if and only if $\Aut(V)$ is generated by automorphisms of $\a^1$-fibrations on $V$.
\end{itemize}
\end{prop}

To determine the structure of $\F_V$ we must have knowledge of \emph{all} possible $1$-standard completions of $V$. Since a general Gizatullin surface might have an immense number of pairwise non-isomorphic $1$-standard completions, there is less hope to compute $\F_V$ in the general case. However, in the following we restrict to surfaces which admit only \emph{rigid} extended divisors. Omiting the precise notion of rigidity of an extended divisor, it means, roughly speaking, that no feather can jump to another boundary component (see \cite{FKZ3} for the notion of rigidity of an extended divisor). For example, if a Gizatullin surface $V$ admits a $1$-standard completion $(X, D)$ such that every feather is attached to an \emph{inner} component, the divisor $D_{\text{ext}}$ is rigid. In this case it follows by the Matching Principle and Lemma \ref{StarComponentAfterReversion} that \emph{any other} extended divisor of $V$ has the same property. Thus, for a surface with this property the dual graph $\Gamma_{D_{\text{ext}}}$ of its extended divisor $D_{\text{ext}}$ can attain at most two different forms. These are precisely the surfaces we study in the following.

Let $V$ be a smooth Gizatullin surface and $(X, D)$ a $1$-standard completion of $V$. Denoting as always by $r_i$ the number of feathers attached to the boundary component $C_i$, we assume for the rest of the article that $V$ has the following property:

\begin{eqnarray*}
&(*)& V \ \text{is a smooth Gizatullin surface such that there exists a $1$-standard completion} \ (X, D) \\
&& \text{of V} \ \text{such that there are no feathers attached to outer components of D and such that} \\
&& r_i > 0 \ \text{holds for at most two} \ i, \ \text{say} \ i = s \ \text{and} \ i = t \ \text{with} \ s \leq t. 
\end{eqnarray*}

\vspace{10pt}
In particular, $(*)$ implies that $r \leq 4$. By the way, if condition $(*)$ holds, the same condition holds for any $m$-standard completion of $V$, since we can change the weight $C^2_1$ of $C_1$ by outer elementary transformations on $C_0$.

To perform explicit calculations with $1$-standard pairs we need to introduce affine coordinates on their completions. One possible way to do this is just to adopt the algorithm in \ref{PresentationGizatullinSurface} for $1$-standard pairs. Thus we can introduce similar affine coordinates $(w_i, z_i)$ and $(x_i, y_i)$, $2 \leq i \leq n$, with the same algorithm, but this time starting with the coordinate system $(u_0, v_0)$ instead of $(x_0, y_0)$. Note, that in this case the map $(x_i, y_i) \mapsto x_i$ describes the correspondence fibration for the pair $(C^\vee_{i^\vee}, C_i)$ if and only if the reversion we deal with is centered in $\lambda = 0$ as well as its inverse. However, since we are interested only in the behaviour of the points blown up in $X \to \f_1$ under birational maps, we don't need the correspondence fibration. Thus such coordinates are sufficient for our purpose, disregarding this small defect.  

In the following we denote by $Q_i$ the i-th configuration invariant of $V$ (\ie, $Q_i$ is the set of base points of feathers attached to $C_i$ modulo the $\c^*$-action on $C_i \backslash (C_{i - 1} \cup C_{i + 1}) \cong \c^*$). Moreover, fixing a (semi-)standard completion  $(X, D)$ of $V$ we denote by $\omega_i := C_i^2$, $2 \leq i \leq n$, the self-intersection numbers of the boundary components.

\begin{prop}\label{TwoFibrations} Let $V$ be a smooth Gizatullin surface satisfying $(*)$ and let $(X, D)$ be a $1$-standard completion of $V$. Then $V$ admits at most two conjugacy classes of $\a^1$-fibrations. Moreover, $V$ admits a unique conjugacy class of $\a^1$-fibrations if and only if $\omega_i = \omega_{i^\vee}$, $t = s^\vee$, $r_s = r_t$ and $Q_s = Q_t$.
\end{prop}

\begin{bew} Since any birational map between $1$-standard completions of $V$ can be decomposed into fibered modifications and reversions and since fibered modifications do not change the invariants $\omega_i$, $r_i$ and $Q_i$, we have to show the following: if $h = \varphi_2 \circ \psi \circ \varphi_1: (X_1, D_1) \dasharrow (X_2, D_2)$, where $\varphi_i$ are reversions and $\psi$ is a fibered modification, then $(X_1 \backslash D_1, \pi_1) \cong (X_2 \backslash D_2, \pi_2)$ as $\a^1$-fibered surfaces. Let us denote by $r^{(1)}_i$ respectively $r^{(2)}_i$ the number of feathers attached to the $i$-th boundary component of $(X_1, D_1)$ respectively $(X_2, D_2)$ (and similarly for the configuration invariants, the boundary components etc.). Since fibered maps do not change the configuration invariants, the Matching Principle yields that $\omega^{(1)}_i = \omega^{(2)}_i$, $r^{(1)}_i = r^{(2)}_i$ and $Q^{(1)}_i = Q^{(2)}_i$ for $i = s, t$.

As we have seen in the proof of Theorem \ref{MainTheorem}, a torus element $t \in \t$, $t(u_0, v_0) = (au_0, bv_0)$ defines a map $(X_1, D_1) \stackrel{\sim}{\to} (X', D')$ such that its restriction on boundary components is of the form

$$h: C^{(1)}_i \backslash C^{(1)}_{i - 1} \stackrel{\sim}{\to} C'_i \backslash C'_{i - 1}, \quad w^{(1)}_i \mapsto w'_i = a^{p_i}b^{q_i}w^{(1)}_i, \quad p_i, q_i, \in \z \quad \text{and} \quad p_iq_j - p_jq_i \neq 0 \ \text{for} \ i \neq j.$$ 

Moreover, elementary shifts $h(u_0, v_0) = (u_0 + av_0^{t - 2}, v_0)$ induce (for general $a$) non-trivial translations $C^{(1)}_{k_t} \backslash C^{(1)}_{k_{t - 1}} \stackrel{\sim}{\to} C^{(2)}_{k_t} \backslash C^{(2)}_{k_{t - 1}}$ in the corresponding $w_i$-coordinates and the identity on all inner components.

Now we proceed as follows. Denote by $A^{(k)}_i$, $k = 1, 2$, the base point sets of the feathers of $(D_k)_{\text{ext}}$ attached to $C^{(k)}_i$. Since $Q^{(1)}_i = Q^{(2)}_i$ for $i = s, t$, we have $A^{(1)}_s = \alpha \cdot A^{(2)}_s$ and $A^{(1)}_t = \beta \cdot A^{(2)}_t$ for some $\alpha, \beta \in \c^*$. Since $p_sq_t - p_tq_s \neq 0$, there exists an element $t = (a, b) \in \t$ such that $a^{p_s}b^{q_s} = \alpha$ and $a^{p_t}b^{q_t} = \beta$. We apply now $t$ to $(X_1, D_1)$ to obtain a new completion $(X', D') = t.(X_1, D_1)$. Let us denote the base point of the component $C^{(2)}_{k_t}$ by $c^{(2)}_{k_t}$ (that is, the point where $C^{(2)}_{k_t}$ is born by an outer blow-up on $C^{(2)}_{k_{t - 1}}$), and similarly we define $c'_{k_t}$. We choose appropriate shifts $h_t(u_0, v_0) := (u_0 + a_iv_0^{t - 2}, v_0)$ which map $c'_{k_t}$ onto $c^{(2)}_{k_t}$ for $t = 2, \dots, r$. Now, the map

$$h := h_r \circ \cdots \circ h_2 \circ t$$

\noindent obviously gives a fibered map $h: (X_1, D_1) \dasharrow (X_2, D_2)$, since all points blown up by $X_1 \to Q$ are mapped onto those blown up by $X_2 \to Q$. It follows that $(X_1 \backslash D_1, \pi_1) \cong (X_2 \backslash D_2, \pi_2)$.

Let us prove the second part. Clearly, if one of the conditions $\omega_i = \omega_{i^\vee}$, $t = s^\vee$, $r_s = r_t$ and $Q_s = Q_t$ does not hold, the $\a^1$-fibered surfaces $(V, \pi)$ and $(V, \pi^\vee)$ cannot be isomorphic. Conversely, let $(X, D)$ be any standard completion of $V$ and let $(X^\vee, D^\vee)$ be the reversed completion (with respect to some center $\lambda \in C^{(1)}_0  \backslash C^{(1)}_1$). Using the arguments above, the condition $\omega_i = \omega_{i^\vee}$, $t = s^\vee$, $r_s = r_t$ and $Q_s = Q_t$ is equivalent to the condition that there is a fibered map which maps the points blown up by $X \to Q$ onto those blown up by $X^\vee \to Q$. Thus the $\a^1$-fibered surfaces $(X \backslash D, \pi)$ and $(X^\vee \backslash D^\vee, \pi^\vee)$ are isomorphic and the corresponding $\a^1$-fibrations are conjugated. It follows that after applying arbitrary fibered maps and reversions to any standard completion $(X, D)$ of $V$ always gives the same conjugacy class of $\a^1$-fibrations. Hence, by Prop. \ref{FactorizationOfBirationalMaps} and Lemma \ref{JonquieresAutomorphisms2}, all $\a^1$-fibrations of $V$ are conjugated.
\end{bew}

\begin{cor}\label{StructureOfAut} Let $V$ be a smooth Gizatullin surface satisfying $(*)$. Then the following holds:
\begin{itemize}
\item[(1)] If $C_3, \dots, C_{n - 1}$ are inner components, the graph $\F_V$ has one of the following structures:

$$\F_V: \begin{xy}
  \xymatrix{
  \bullet \ar@{<->}[r] & \bullet \\
  }
\end{xy} \quad \text{or} \quad \F_V: \bullet \rcirclearrowleft .$$

\item[(2)] If one of the components $C_3, \dots, C_{n - 1}$ is an outer component, the graph $\F_V$ has one of the following structures:
\begin{itemize}
\item[(i)] $\F_V$ consists of two vertices $v_1$ and $v_2$ and uncountable many arrows between $v_1$ and $v_2$.
\item[(ii)] $\F_V$ consists of a unique vertex $v$ and uncountable many arrows starting and ending at $v$.
\end{itemize}
\end{itemize}
\noindent In both cases $\F_V$ admits a unique vertex if and only if $\omega_i = \omega_{i^\vee}$, $t = s^\vee$, $r_s = r_t$ and $Q_s = Q_t$.
\end{cor}

\begin{bew} The last assertion follows directly from Proposition \ref{TwoFibrations} and the fact that the vertices of $\F_V$ are in bijection to conjugacy classes of $\a^1$-fibrations of $V$.
 
Assertion (1) is proven in \cite{Ko}, Theorem 3.4. Let us show (2). By Proposition \ref{TwoFibrations} it follows that $\F_V$ has at most two vertices. We show that $\F_V$ has uncountable many arrows. Let $(X, D)$ be an arbitrary $1$-standard completion of $V$. It is sufficient to show that the group $\Aut(X, D)$ is finite. Indeed, two reversions $(X, D) \dasharrow (X', D')$ and $(X, D) \dasharrow (X'', D'')$ centred in $\lambda_1$ and $\lambda_2$ define the same arrow if and only if there is an automorphism $\varphi \in \Aut(X, D)$ which maps $\lambda_1$ to $\lambda_2$. Therefore there are uncountable many pairwise non-equivalent reversions starting from $(X, D)$.

Elements $\varphi \in \Aut(X, D)$ are lifts of maps of the form 

$$\varphi(u_0, v_0) = (au_0 + \alpha v_0 + \beta, bv_0 + c).$$

\noindent Since $\varphi$ stabilizes $C_2$ (which is the closure of $\{ v_0 = 0 \}$) we have $c = 0$. The map $\varphi$ induces on the two (inner) boundary components $C_s$ and $C_t$ maps of the form

\begin{eqnarray*}
&& C_s \backslash (C_{s - 1} \cup C_{s + 1}) \to C_s \backslash (C_{s - 1} \cup C_{s + 1}), \quad w_s \mapsto a^{p_s}b^{q_s}w_s, \\
&& C_t \backslash (C_{t - 1} \cup C_{t + 1}) \to C_t \backslash (C_{t - 1} \cup C_{t + 1}), \quad w_t \mapsto a^{p_t}b^{q_t}w_t, \\
\end{eqnarray*}

\noindent such that $p_sq_t - p_tq_s \neq 0$. Thus, since $\varphi$ stabilizes the base point sets of the feathers, $a$ and $b$ can only attain a finite number of values. Moreover, since $\varphi$ stabilizes $C_2 \cap C_3$, we have $(ac_{k_3} + \beta, 0) = (c_{k_3}, 0)$. Thus $\beta$ is determined by the value of $a$. Finally, we consider the blow-up $X_2 \to X_1$ from \ref{PresentationGizatullinSurface} (2), which creates the outer component $C_{k_3}$ by a blow-up in $(u_0, v_0) = (c_{k_3}, 0)$. Introducing the coordinates 

$$(w_{k_3}, z_{k_3}) = \left( \frac{u_0 - c_{k_3}}{v_0}, v_0 \right)$$ 

\noindent gives $C_{k_3} \backslash C_{k_3 - 1} = \{ z_{k_3} = 0 \}$ and $\varphi$ induces the linear map

$$C_{k_3} \backslash C_{k_3 - 1} \stackrel{\sim}{\to} C_{k_3} \backslash C_{k_3 - 1}, \quad w_{k_3} \mapsto ab^{-1}w_{k_3} + \alpha b^{-1}.$$

\noindent Since $\varphi$ stabilizes the intersection point $C_{k_3} \cap C_{k_3 + 1}$, which corresponds to a certain point $(w_{k_3}, z_{k_3}) = (c'_{k_4}, 0)$, we have $ab^{-1}c'_{k_4} + \alpha b^{-1} = c'_{k_4}$ and thus $\alpha$ is determined by the values of $a$ and $b$. It follows that $\Aut(X, D)$ is finite, which completes the proof.
\end{bew}

Assume that $V$ is a Gizatullin surface such that the boundary divisor $D$ in standard form admits a component $C_i$, $i \geq 2$, with self-intersection $\leq -3$. According to Prop. \ref{SurfaceTree} we have a natural exact sequence

$$0 \to H \to \Aut(V) \to \Pi_1(\F_V) \to 0.$$

By \cite{BD2}, Prop. 2.7 and Remark 2.8 it follows that the image of any algebraic subgroup of $\Aut(V)$ in $\Pi_1(\F_V)$ is finite. More precise, there is at most one nontrivial element in the image, which consists, if it exists, of one path of the form $\varphi^{-1}\sigma\varphi$, where $\varphi$ is a path from a vertex $[(X, D)]$ to another vertex $[(X', D')]$ and $\sigma$ is a loop of length $1$ based at the vertex $[(X', D')]$ and representing a reversion $(X', D') \dasharrow (X'', D'')$ between two isomorphic pairs. The last case is only possible if $C^2_i = C^2_{i^ \vee}$, otherwise the image of every algebraic supgroup of $\Aut(V)$ in $\Pi_1(\F_V)$ is trivial.

In particular, if $\Pi_1(\F_V)$ is not countable, then $\Aut(V)$ is not generated by a countable set of algebraic subgroups.

Since any Gizatullin surface as in Corollary \ref{StructureOfAut} (2) admits a boundary component with self-intersection $\leq -3$, we immediately obtain the following corollary:

\begin{cor}\label{NotCountableGenerated} Let $V$ be a Gizatullin surface as in Corollary \ref{StructureOfAut} (2). Then $\Aut(V)$ is not generated by a countable set of algebraic subgroups.
\end{cor}

We assume further that a smooth Gizatullin surface $V$ satisfies $(*)$ with $r = 4$. In other words, a $1$-standard completion $(X, D)$ of $V$ admits precisely three outer components in $D^{\geq 2}$, namely $C_{k_2} = C_2, C_{k_3}$ and $C_{k_4} = C_n$, and precisely two families of feathers, both attached to inner components of $D$. Though this consition is not really necessary (\ie the result remains valid without the assumption $r = 4$), it simplifies the proof a litte. Furthermore, we assume that $\F_V$ consists of two vertices. Under these special assumptions we have

\begin{cor}\label{UncountableFreeGroup} Denote by $\Aut(V)_{\text{alg}}$ the (normal) subgroup of $\Aut(V)$ generated by all algebraic subgroups of $\Aut(V)$. Then $\Aut(V)/\Aut(V)_{\text{alg}}$ contains a free group $F$ over an uncountable set of generators.
\end{cor}

\begin{bew} First we show the following\\

\noindent \textbf{Claim:} Let $(X_1, D_1)$ and $(X_2, D_2)$ two $1$-standard completions of $V$ such that $(X_1 \backslash D_1, \pi_1) \cong (X_2 \backslash D_2, \pi_2)$. Then $(X_1, D_1) \cong (X_2, D_2)$.\\

We consider any isomorphism $\varphi: (X_1 \backslash D_1, \pi_1) \cong (X_2 \backslash D_2, \pi_2)$, which is by Lemma \ref{JonquieresAutomorphisms2} a lift of some automorphism of the form $\varphi(u_0, v_0) = (au_0 + P(v_0), bv_0 + c)$. Let us write $P(v_0) = \sum_{j = 0}^m a_jv_0^j$ with $m \geq 1$ and $a_j \in \c$. Since the term $a_jv_0^j$ induces only motions on outer components $C_{k_i}$ with $i \geq j + 2$, it follows that $\tilde{\varphi}(u_0, v_0) := (au_0 + a_0 + a_1v_0, bv_0 + c)$ also gives a fibered map $(X_1 \backslash D_1, \pi_1) \stackrel{\sim}{\to} (X_2 \backslash D_2, \pi_2)$, and furthermore, \emph{all} fibered maps $(X_1 \backslash D_1, \pi_1) \stackrel{\sim}{\to} (X_2 \backslash D_2, \pi_2)$ are given in this way (this is where we need $r = 4$). But again by Lemma \ref{JonquieresAutomorphisms2} such maps already define isomorphisms $(X_1, D_1) \stackrel{\sim}{\to} (X_2, D_2)$, since $\tilde{\varphi} \in \Jon \cap \Aff$.\\

Roughly speaking, the idea behind the following is that we subdivide the arrows of $\F_V$ into (an uncountable set of) disjoint pairs. Every such pair defines a loop in $\F_V$ and we associate to every such loop an automorphism of $V$. The group generated by these automorphism will be our $F$.

We choose an uncountable subset $A \subseteq \c$ as follows: By the above claim there exists a unique $1$-standard completion of $V$ representing a given vertex. Choose two $1$-standard completions of $V$ representing the vertices of $\F_V$, say $(X, D)$ and $(X', D')$. By the proof of Proposition \ref{StructureOfAut} the group $G := \Aut(X, D)$ is finite. Identifying $C_0 \backslash C_1 \cong \c$, we choose $\tilde{A} \subseteq \c$ to be the set of representatives of the equivalence classes of $\c$ modulo the $G$-action (in other words, $\tilde{A}$ contains precisely one element of each $G$-orbit of $\c$). Then $\tilde{A}$ parametrizes the arrows in $\F_V$. Similarly we introduce a subset $\tilde{A}' \subseteq C'_0 \backslash C'_1 \cong \c$, which also parametrizes the arrows in $\F_V$, since every arrow in $\F_V$ has an inverse. Thus there exists a bijection $f: \tilde{A} \stackrel{\sim}{\to} \tilde{A}'$ which satisfies the follwoing property: if a reversion $\psi: (X, D) \dasharrow (X', D')$ is centred in $a \in \tilde{A}$, its inverse $\psi^{-1}\: (X', D') \dasharrow (X, D)$ is centred in $f(a) \in \tilde{A}'$. Now choose an arbitrary involution $\sigma: \tilde{A} \to \tilde{A}$ such that $\sigma(a) \neq a$ for all $a \in \tilde{A}$. By identifying $\tilde{A}$ with $\tilde{A}'$ via $f$ we obtain an involution $\sigma': \tilde{A}' \to \tilde{A}'$ (with $\sigma'(a') \neq a'$ for all $a' \in \tilde{A}'$) such that $f \circ \sigma = \sigma' \circ f$. We let $A \subseteq \tilde{A}$ be a set of representatives of the equivalence relation $\sim_\sigma$ on $\tilde{A}$, where $a \sim_\sigma b$ holds if $a = b$ or $b = \sigma(a)$. Furthermore, we let $A' := f(A) \subseteq \tilde{A}'$ (which is a set of representatives of the equivalence relation $\sim_{\sigma'}$ on $\tilde{A}'$ induces by $\sigma'$).

For $a, a' \in \c$ we denote by $\psi_a: (X, D) \dasharrow (X', D')$ a reversion centred in $a$ and by $\psi'_{a'}: (X', D') \dasharrow (X, D)$ a reversion centred in $a'$. Now, for any $a \in A$ we consider the birational map

$$\varphi_a := \psi'_{(\sigma' \circ f)(a)} \circ \psi_a: (X, D) \dasharrow (X', D') \dasharrow (X, D).$$

\noindent By construction, its inverse $\varphi^{-1}_a := \psi'_{f(a)} \circ \psi_{\sigma(a)}$ is centred in $\sigma(a)$. The map $\varphi_a$ restricts to an automorphism of $V$, which will be denoted by the same letter. We claim that the group $F := \langle \varphi_a \mid a \in A \rangle$ is free and intersects $\Aut(V)_{\text{alg}}$ trivially.

First we show the freeness. By the uniqueness part of Prop. \ref{FactorizationOfBirationalMaps}, any composition of the form

$$(\varphi_{a_1})^{\delta_1} \circ \cdots \circ (\varphi_{a_s})^{\delta_s}: (X, D) \dasharrow (X, D) \quad \text{with} \quad \delta_i \in \z \backslash \{ 0 \} \quad \text{and} \quad a_i \neq a_{i + 1}$$

\noindent has minimal length. This implies that any automorphism $\varphi := (\varphi_{a_1}\vert_V)^{\delta_1} \circ \cdots \circ (\varphi_{a_s}\vert_V)^{\delta_s} \in \Aut(V)$ with $s \geq 1$, $\delta_i \in \z \backslash \{ 0 \}$ and $a_i \neq a_{i + 1}$ is not trivial.

Now, the image of $\varphi$ in $\Pi_1(\F_V)$ consists of a product of loops based on $[(X, D)]$ of length $\geq 2$. Since the image of any element of $\Aut(V)_{\text{alg}}$ can only contain loops of length $1$, we obtain that $F \cap \Aut(V)_{\text{alg}} = \{ \id_V \}$ and the proof is complete.
\end{bew}

Let us give a concrete example of a smooth non-homogeneous Gizatullin surface with a huge automorphism group:

\begin{ex} A smooth Gizatullin surface $V$ with the following extended divisor is non-homogeneous and has a huge automorphism group:

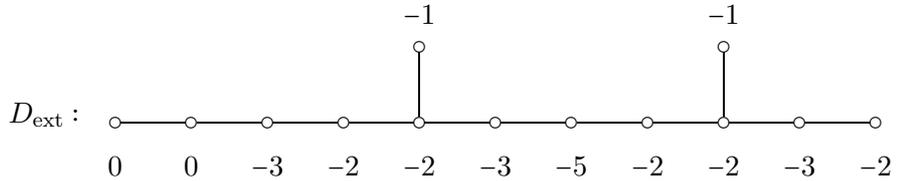
\vspace{15pt}
$$D_{\text{ext}}: \quad \cu{0} \lin \cu{0} \lin \cu{-3} \lin \cu{-2} \lin \cu{-2} \nlin \cshiftup{-1}{} \lin \cu{-3} \lin \cu{-5} \lin \cu{-2} \lin \cu{-2} \nlin \cshiftup{-1}{} \lin \cu{-3} \lin \cu{-2} \quad .$$

\noindent Indeed, it is easy to check that the only outer components are $C_2, C_6$ and $C_{10}$, and moreover, $\mathfrak{E}_D = \{ 5, 9 \}$ and $\mathfrak{E}_{D^\vee} = \{ 3, 5, 7, 9 \}$, hence $\mathfrak{E}^\vee_{D^\vee} = \{ 3, 5, 7, 9 \}$. Thus by Theorem \ref{MainTheorem}, $\Aut(V)$ admits at least two fixed points (the intersection points of the feathers with their duals) and Cor. \ref{NotCountableGenerated} and Cor. \ref{UncountableFreeGroup} give the desired statement about the hugeness of $\Aut(V)$. 
\end{ex}

\normalsize\rmfamily

\noindent Albert-Ludwigs-Universit\"at Freiburg, Mathematisches Institut, Eckerstrasse 1, 79104 Freiburg, Germany\\
\noindent \textit{E-mail address:} sergei.kovalenko@math.uni-freiburg.de 
 
\end{sloppypar}


\begin{thebibliography}{KMMH}
\bibitem[BD1]{BD1} J. Blanc, A. Dubouloz, \textit{Automorphisms of $\a^1$-fibered affine surfaces}, Trans. Amer. Math. Soc. 363 (2011), no. 11, 5887-5924.
\bibitem[BD2]{BD2} J. Blanc, A. Dubouloz, \textit{Affine surfaces with a huge group of automorphisms}, Int. Math. Res. Not. 2013, doi: 10.1093/imrn/rnt202
\bibitem[Be]{Be} J. Bertin, \textit{Pinceaux de droites et automorphismes des surfaces affines}, J. Reine Angew. Math. 341 (1983), 32 - 53.
\bibitem[BML]{BML} T. Bandman, L. Makar-Limanov, \textit{Affine surfaces with $AK(S) = \c$}, Michigan Math. J. 49 (2001), 567 - 582. 
\bibitem[Da]{Da} D. Daigle, \textit{Classification of linear weighted graphs up to blowing-up and blowing-down}, Canad. J. Math. 60 (2008), no. 1, 64-87.
\bibitem[DG1]{DG1} V. I. Danilov, M. H. Gizatullin, \textit{Examples of nonhomogeneous quasihomogeneous surfaces}, Mathematics of the USSR-Izvestiya 8:1 (1974), 43-60.
\bibitem[DG2]{DG2} V. I. Danilov, M. H. Gizatullin, \textit{Automorphisms of affine surfaces I}, Izv. Akad. Nauk SSSR Ser. Mat. 39:3 (1975), 523-565.
\bibitem[DG3]{DG3} V. I. Danilov, M. H. Gizatullin, \textit{Automorphisms of affine surfaces II}, Izv. Akad. Nauk SSSR Ser. Mat. 41:1 (1977), 54-103.
\bibitem[Du]{Du} A. Dubouloz, \textit{Completions of normal affine surfaces  with a trivial Makar-Limanov invariant}, Michigan Math. J. 52 (2004), 289 - 308.
\bibitem[Fr]{Fr} G. Freudenburg, \textit{Algebraic theory of locally nilpotent derivations}, Encyclopaedia of Mathematical Sciences 136, Invariant Theory and Algebraic Transformation Groups 7, Springer (2006).
\bibitem[FZ1]{FZ1} H. Flenner, M. Zaidenberg, \textit{Normal affine surfaces with $\c^*$-actions}, Osaka J. Math. 40 (2003), no. 4, 981-1009.
\bibitem[FZ2]{FZ2} H. Flenner, M. Zaidenberg, \textit{Locally nilpotent derivations on affine surfaces with a $\c^*$-actions}, Osaka J. Nath. 42 (2005), no. 4, 931-974.
\bibitem[FKZ1]{FKZ1} H. Flenner, S. Kaliman, M. Zaidenberg, \textit{Birational transformations of weighted graphs}, Affine algebraic geometry, 107-147, Osaka Univ. Press, Osaka, 2007.
\bibitem[FKZ1C]{FKZ1C} H. Flenner, S. Kaliman, M. Zaidenberg, \textit{Corrigendum to Our Paper "Birational transformations of weighted graphs"}, Affine Algebraic Geometry, 35-38, CRM Proc. Lecture Notes, 54, Amer. Math. Soc., Providence, RI, 2011.
\bibitem[FKZ2]{FKZ2} H. Flenner, S. Kaliman, M. Zaidenberg, \textit{Completions of $\c^*$-surfaces}, Affine algebraic geometry, 149-201, Osaka Univ. Press, Osaka, 2007.
\bibitem[FKZ3]{FKZ3} H. Flenner, S. Kaliman, M. Zaidenberg, \textit{Uniqueness of $\c^*$- and $\c_+$-actions on Gizatullin surfaces}, Transform. Groups 13 (2008), no. 2, 305-354.
\bibitem[FKZ4]{FKZ4} H. Flenner, S. Kaliman, M. Zaidenberg, \textit{Smooth affine surfaces with non-unique $\c^*$-actions}, J. Algebraic Geom. 20 (2011), no. 2, 329-398.
\bibitem[Gi]{Gi} M. H. Gizatullin, \textit{Quasihomogeneous affine surfaces}, Izv. Akad. Nauk SSSR Ser. Mat. 35 (1971), 1047 - 1071.
\bibitem[Ha]{Ha} R. Hartshorne, \textit{Algebraic Geometry}, Graduate texts in mathematics: 52, Springer (2006).  
\bibitem[Ko]{Ko} S. Kovalenko, \textit{Transitivity of automorphism groups of Gizatullin surfaces}, arXiv: 1304.7116.
\bibitem[Mi]{Mi} M. Miyanishi, \textit{Open algebraic surfaces}. CRM Monograph Series, 12. American Mathematical Society, Providence, RI (2001).
\bibitem[Se]{Se} J.-P- Serre, \textit{Trees}. Springer-Verlag Berlin Heidelberg New York 1980.
\end{thebibliography}
\end{document}